\documentclass[reqno,12pt]{amsart}

\usepackage[initials]{amsrefs}
\usepackage{amscd,amssymb,comment,euscript}
\usepackage{graphicx}
\usepackage{enumerate}

\usepackage[all]{xy}

\theoremstyle{definition}
\newtheorem{notation}{Notation}[section]

\theoremstyle{plain}

\newtheorem{prop}[notation]{Proposition}
\newtheorem{lemma}[notation]{Lemma}
\newtheorem{thm}[notation]{Theorem}
\newtheorem{coro}[notation]{Corollary}

\theoremstyle{definition}
\newtheorem{definition}[notation]{Definition}

\theoremstyle{definition}
\newtheorem{remark}[notation]{Remark}

\newcommand{\SA}{\textnormal{$*$-SubAlg}}
\newcommand{\Sb}{\textnormal{Stab}}
\newcommand{\Aut}{\operatorname{Aut}}
\newcommand{\Ad}{\operatorname{Ad}}
\newcommand{\Inn}{\operatorname{Inn}}
\newcommand{\lcm}{\operatorname{lcm}}
\newcommand{\dist}{\operatorname{dist}}
\newcommand{\dH}{\operatorname{d_H}}
\newcommand\id{{\operatorname{id}}}

\begin{document}

\title[Primitivity of free products]{Primitivity of unital full free products of residually finite dimensional $C^*$-algebras}

\author[Dykema]{Ken Dykema}
\author[Torres--Ayala]{Francisco Torres-Ayala}
\address{Department of Mathematics, Texas A\&M University,
College Station, TX 77843-3368, USA}
\email{kdykema@math.tamu.edu}
\email{francisc@math.tamu.edu}
\thanks{Research supported in part by NSF grant DMS--0901220.}

\date{June 19, 2012}

\subjclass[2000]{46L09 (46L05)}

\keywords{Primitive C$^*$--algebra, Full free product}

\begin{abstract}
A $C^*$-algebra is called primitive if it admits a faithful and irreducible $*$--representation. 
We show that if $A_1$ and $A_2$ are separable, unital, residually finite dimensional C*-algebras satisfying
$(\dim(A_1)-1)(\dim(A_2)-1) \geq 2$, then the unital $C^*$-algebra full free product,  $A=A_1*A_2$,
is primitive. It follows  that $A$ is  antiliminal, it has an uncountable
family of pairwise inequivalent irreducible faithful $*$--representations   and  the set of 
pure states is  w*-dense in the state space.
\end{abstract}

\maketitle

\section{Introduction}

A $C^*$-algebra is called primitive if it admits a  faithful and irreducible $*$--representation. Thus the simplest examples are matrix algebras. A nontrivial example, shown independently by Choi and Yoshizawa,  is the full group $C^*$-algebra of the free group on $n$ elements, $2\leq n\leq \infty$, see \cite{Choi} and \cite{Yoshizawa}.
In \cite{Murphy}, Murphy gave numerous conditions for primitivity of full group $C^*$-algebras.
More recently, T. \AA. Omland showed in  ~\cite{Omland-PrimitivityConditionsTwistedGroups} that for
$G_1$ and $G_2$ countable amenable discrete groups and $\sigma$ a multiplier
on the free product $G_1*G_2$, the full twisted
group C*-algebra $C^*(G_1*G_2,\sigma)$ is primitive whenever $(|G_1|-1)(|G_2|-1)\geq 2$.

We prove that given two nontrivial,  separable, unital, residually finite dimensional  $C^*$-algebras $A_1$ and $A_2$,  their  unital $C^*$-algebra full free product $A_1*A_2$ is primitive except when $A_1=\mathbb{C}^2=A_2$.
The methods  used are essentially different from those in ~\cite{Murphy}, ~\cite{Bedos&Omland-Modular}, ~\cite{Bedos&Omland-Amenable} and ~\cite{Omland-PrimitivityConditionsTwistedGroups} but do rely on Exel and Loring's result~\cite{Exel&Loring} that $A_1*A_2$
is itself residually finite dimensional.
Roughly speaking, we first show that if $(\dim(A_1)-1)(\dim(A_2)-1)\geq 2$, then there is an abundance of irreducible finite dimensional $*$--representations and later, by means of a sequence of  approximations, we construct an irreducible and faithful $*$--representation.

The paper is divided as follows.  Section \ref{Automorphisms} recalls some  facts about $*$--automorphisms of  finite dimensional $C^*$-algebras. Section \ref{LieGroups} recalls some known result on Lie groups that will be used later. Section \ref{Intersection&Perturbation} is fully  devoted in  proving  Theorem \ref{DensityOfSmallIntersection} which is about perturbing a pair of proper unital $C^*$-subalgebras of a matrix algebra in such a way that they have trivial intersection.  Theorem \ref{DensityOfSmallIntersection}  is the cornerstone for the rest of the results in the paper.  Lastly, section \ref{PrimitivitySection} contains the proof of the main theorem about primitivity and  some  consequences.

\begin{notation}
Given a Hilbert space $H$, we denote the set of bounded linear operators by $\mathbb{B}(H)$ and the set of compact operators by $\mathbb{K}(H)$.

 For a unital $C^*$-algebra $A$, $\SA(A)$ denotes the set of all  unital $C^*$-subalgebras of $A$ and $\mathbb{U}(A)$ denotes the set of unitary elements of $A$. For simplicity, given a Hilbert space $H$ we write $\mathbb{U}(H)$ instead of $\mathbb{U}(\mathbb{B}(H))$.

 By $\Aut(A)$ we denote the set of $*$--automorphisms of $A$. For $u$ in $\mathbb{U}(A)$  we let $\Ad u$ denote the $*$--automorphism of $A$ given by $\Ad u(x)=uxu^*$. The set of all $*$--automorphisms of the form $\Ad u$, for some $u$, is called the set of inner automorphism and it is denoted by $\Inn(A)$.

For a unital $C^*$-algebra $A$, $C(A)$ denotes its center. In other words
$$
C(A)=\{x\in A: xa=ax \quad \textrm{for all $a\in A$} \}.
$$

For a positive integer $n$, $M_n$ denotes the set of $n\times n$ matrices over $\mathbb{C}$
and $S_n$ denotes the  permutation group of the set $\{1,\dots, n\}$.

\end{notation}

\section{$*$-Automorphisms of finite dimensional $C^*$-algebras}
\label{Automorphisms}

By a $*$--automorphism of a $C^*$-algebra we mean a  bijective map,  from the algebra onto itself, that is linear and preserves products and  adjoints. 

In this section we recall some basic results concerning $*$--auto\-morph\-isms of finite dimensional $C^*$-algebras and in particular a precise  algebraic relation between the group of  $*$--automorphism and the subgroup of inner $*$--automorphisms.

Any $*$-homomorphism from a simple $C^*$-algebra is either zero or  injective (since its kernel is an ideal). Even more,  any non-zero $*$-endomorphism of a finite dimensional simple $C^*$-algebra is a  $*$--auto\-morph\-ism. Indeed, any such $*$-endomorphism is injective and thus it is bijective (by finite  dimensionality) and a straightforward computation shows its inverse
is a $*$-endomorphism.
As a consequence  any  $*$--automorphism of a finite dimensional $C^*$-algebra moves, without breaking, each  one of its simple $C^*$-subalgebras (we may think this as blocks) with the same dimension. Thus modulo an inner $*$--automorphism,  a $*$--automorphism is  just a permutation. We make the last statement precise with the following two propositions.

\begin{prop}\label{2ndDecompositionAutomorphisms}
Let $B$ be a finite dimensional $C^*$-algebra and assume $B$ decomposes as
$$
\oplus_{j=1}^J B_j
$$
and there is a positive integer $n$ such that  all $B_j$ are $*$--isomorphic to $M_n$.

Fix $\{ \beta_j: B_j \to M_n \}_{1\leq j \leq J}$ a set of  $*$--isomorphisms. 

\begin{enumerate}
\item For a permutation $\sigma$ in $S_J$ define $\psi_{\sigma} : B \to B$ by
$$
\psi_{\sigma}(b_1,\dots , b_J) = ( \beta_1^{-1} \circ \beta_{\sigma^{-1}(1)} (b_{\sigma^{-1}(1)}), \dots ,\beta_J^{-1} \circ \beta_{\sigma^{-1}(J)} (b_{\sigma^{-1}(J)}))
$$
Then $\psi_{\sigma}$ lies  in $\Aut(B)$ and the map $\sigma \mapsto \psi_{\sigma}$ defines a group  embedding of $S_J$ into $\Aut(B)$.

\item Every element  $\alpha$ in $\Aut(B)$ factors as
$$
\left( \oplus_{j=1}^J \Ad u_j \right) \circ \psi_{\sigma}
$$
for some permutation $\sigma$ in $S_J$ and unitaries $u_j $ in $\mathbb{U}(B_j)$.

\item There is a exact sequence 
$$
0 \to \Inn(B) \to \Aut(B) \to S_J\to 0.
$$

\end{enumerate}
\end{prop}

So far we have consider $C^*$-algebras with only one type of block subalgebra, so to speak. Next proposition shows that a $*$--automorphism can not mix  blocks of different dimensions. As a consequence, and  along with Proposition \ref{2ndDecompositionAutomorphisms}, we get a general decomposition of $*$--automorph\-isms of finite dimensional $C^*$-algebras.

\begin{prop}\label{ThirdDecompositionAutomorphisms}
Let $B$ be a finite dimensional $C^*$-algebra. Decompose $B$ as
$$
\oplus_{i=1}^I\oplus_{j=1}^{J_i} B(i,j)
$$
where for each $i$, there is a positive integer  $n_i $ such that $B(i,j)$ is isomorphic to $M_{n_i}$ for all $1\leq j \leq  J_i$, i.e.
we group   subalgebras  that are isomorphic to the same matrix algebra, and where $n_1<n_2<\cdots<n_I$..

Then any   $\alpha$  in $\Aut(B)$  factors as $\alpha = \oplus_{i=1}^I \alpha_i$ where
$$
\alpha_i : \oplus_{j=1}^{J_i} B(i,j) \to \oplus_{j=1}^{J_i} B(i,j)
$$
is a  $*$--isomorphism.  
\end{prop}

\section{Useful results from Lie Groups}\label{LieGroups}

In this section we  summarize some result that, later on,  will be repeatedly used. Definitions and proofs of results mentioned in this section can be found in ~\cite{Kirillov} and ~\cite{Helgason}.

The next two theorems are quite  important and will be used  in the next section.  

\begin{thm}
Any closed subgroup of a Lie group is a Lie subgroup.
\end{thm}

\begin{thm}
 Let $G$ be a Lie group of dimension $n$ and $H \subseteq G$ be  a Lie subgroup of dimension $k$. 
\begin{enumerate} 
\item  Then the  left coset space $G/H$   has a natural structure of a manifold of dimension $n-k$ such that the 
canonical quotient map $\pi : G \to G/H $, is a fiber bundle, with fiber diffeomorphic to $H$. 
\item If $H$ is a normal Lie subgroup then $G/H$ has a canonical structure of a Lie group. 
\end{enumerate}
\end{thm}

The next proposition is from Corollary 2.21 in~\cite{Kirillov}.  

\begin{prop}\label{Orbit&Stabilizer}
Let $G$ denote a Lie group and assume it acts smooth\-ly on  a manifold $M$. For $m\in M$ 
let $\mathcal{O}(m)$ denote its orbit and  $\Sb(m)$ denote its stabilizer i.e.
\begin{eqnarray*}
\mathcal{O}(m) & = & \{ g.m: g\in G \}, \\
\Sb(m) &= &\{  g \in G: g.m=m\}.
\end{eqnarray*}
The orbit $\mathcal{O}(m)$ is an immersed submanifold of $M$. If $\mathcal{O}(m)$ is compact, then the map $g \mapsto g.m $, is a  diffeomorphism from $G/\Sb(m)$ onto $\mathcal{O}(m)$.
(In this case we say $\mathcal{O}(m)$ is an embedded submanifold of $M$.)
\end{prop}

\begin{coro}\label{ProductOfLieSubgroups}
Let $G$ be a compact Lie group and let $K$ and $L$ be  closed subgroups of $G$. The  subspace
$KL=\{ kl : k \in K, l\in L \}$ is an  embedded submanifold of $G$ of dimension
$$
\dim K + \dim L - \dim (L\cap K).
$$
\end{coro}

\begin{proof}
First of all $KL$ is compact. This follows from the fact that multiplication is continuous and both $K$ and $L$ are compact.
Consider the action of $K\times L$ on $G$ given by $(k,l).g=kgl^{-1}. $
Notice that the orbit of $e$ is precisely  $KL$. By Proposition \ref{Orbit&Stabilizer}, $KL$ is an immersed submanifold  diffeomorphic to 
$K\times L / \Sb(e).$
Since it is compact, it is an embedded submanifold.
But
$\Sb(e)=\{ (x,x): x \in K \cap L \}$
and we conclude 
$$
\dim KL = \dim (K \times L) -\dim \Sb(e)=\dim K + \dim L -\dim(K\cap L) .
$$
\end{proof}

\begin{prop}\label{LocalCrossSection}
Let $G$ be a compact Lie group and let $H$ be a closed subgroup. Let $\pi$ denote the quotient map onto $G/H$.

There are:
\begin{enumerate}

\item $\mathcal{N}_G$, a compact neighborhood of $e$ in $G$,
\item $\mathcal{N}_H$, a compact neighborhood of $e$ in $H$,
\item $\mathcal{N}_{G/H}$, a compact neighborhood of $\pi(e)$ in $G/H$,
\item a continuous function $s :\mathcal{N}_{G/H}(\pi(e)) \to G$ satisfying
\begin{enumerate}
\item $s(\pi(e))=e$ and $\pi(s(y))=y$ for all $y $ in $\mathcal{N}_{G/H}(\pi(e))$, 
\item The map
\begin{eqnarray*}
\mathcal{N}_H \times \mathcal{N}_{G/H} & \to & \mathcal{N}_G, \\
(h,y) & \mapsto & hs_g(y)
\end{eqnarray*}
is a homeomorphism.
\end{enumerate}
\end{enumerate}

\end{prop}

\begin{proof}
Let $\mathfrak{g}$ and $\mathfrak{h}$ denote, respectively,  the Lie algebras of $G$ and $H$. Take $\mathfrak{m}$ a  vector subspace such that $\mathfrak{g}$ is the direct sum of $\mathfrak{h}$ and $\mathfrak{m}$. By Lemmas  2.4 and 4.1 in \cite{Helgason}, chapter 2, there are compact neighborhoods  $U_{\mathfrak{g}}$, $U_{\mathfrak{h}}$ and $U_{\mathfrak{m}}$ of $0$ in $\mathfrak{g}$, $\mathfrak{h}$ and $\mathfrak{m}$, respectively, such that  the map
\begin{eqnarray*}
U_{\mathfrak{m} }\times U_{\mathfrak{h}} &\to & U_{\mathfrak{g}},\\
(a,b) & \mapsto & \exp(a)\exp(b)
\end{eqnarray*}
is an homeomorphism and $\pi$ maps homeomorphically $\exp(U_{\mathfrak{m}})$ onto a compact neighborhood of $\pi(e)$. Call the latter neighborhood $\mathcal{N}_{G/H}$.
Take $\mathcal{N}_G=exp(U_{\mathfrak{g}})$, $\mathcal{N}_H=\exp(U_{\mathfrak{h}})$ and $s$ the inverse of $\pi$ restricted to $\exp(U_{\mathfrak{m}})$.

\end{proof}

\section{Intersection of Finite Dimensional $C^*$-algebras and Perturbations }
\label{Intersection&Perturbation}

In this section we fix a positive integer $N$ and, unless stated otherwise, $B_1 \varsubsetneq M_N$ and $B_2 \varsubsetneq M_N$ denote proper unital  $C^*$-subalgebras of $M_N$.

The main purpose if this section is give a proof of the following theorem (recall that for a $C^*$-algebra $A$, $C(A)$ denotes its center).

\begin{thm}\label{DensityOfSmallIntersection}
Assume one of the following conditions holds:
\begin{enumerate}
\item\label{SimpleCaseDensityOfSmallIntersection} $\dim C(B_1)=1=\dim C(B_2)$,
\item \label{BigCenter*SimpleCaseDensityOfSmallIntersection} $\dim C(B_1)  \geq 2$, $\dim C(B_2)=1$ and $B_1$ is $*$--isomorphic to
$$
M_{N/\dim C(B_1)}\oplus \cdots \oplus M_{N/\dim C(B_1)},
$$
\item\label{CenterCenter2CaseDensityOfSmallIntersection} $\dim C(B_1) = 2=\dim C(B_2)$, $B_1$ is $*$--isomorphic to
$$
M_{N/2} \oplus M_{N/2},
$$
 and $B_2$ is $*$--isomorphic to
$$
M_{N/2}\oplus M_{N/(2k)}
$$
where $k \geq 2$,
\item\label{BigCenterCaseDensityOfSmallIntersection}$\dim C(B_1)\geq 2$, $\dim C(B_2) \geq 3$  and, for $i=1,2$, $B_i$ is $*$--isomorphic to
$$
M_{N/\dim C(B_i)}\oplus \cdots \oplus M_{N/ \dim C(B_i)}.
$$
\end{enumerate} 
Then 
$$\Delta(B_1,B_2):=\{  u\in \mathbb{U}(M_N): B_1\cap uB_2u^* = \mathbb{C} \}$$
 is  dense in $\mathbb{U}(M_N)$.
\end{thm}

The $C^*$-algebra $uB_2u^*$ is what we call a perturbation of $B_2$ by $u$. With this nomenclature we are trying to prove that, in the cases mentioned above,  almost  always  we can  perturb one $C^*$-subalgebra a little bit  in such a way that the intersection with the other one is the smallest possible. 
 
Roughly speaking, the idea behind  is to show that the complement of  $\Delta(B_2,B_2)$  can be locally parametrized with strictly fewer variables than $\dim \mathbb{U}(M_N)=N^2$. Thus, the complement of $\Delta(B_1,B_2)$ is, topologically speaking,  small.

We start with some definitions.
The group $\mathbb{U}(B_1)$  acts on $\SA(B_1)$  via 
$(u,B)\mapsto uBu^*$
and the equivalence relation on $\SA(B_1)$ induced by this action will be denoted by $\sim_{B_1}$.  Specifically, we have 
$$
B \sim_{B_1} C  \Leftrightarrow \exists u \in \mathbb{U}(B_1)  : uBu^*=C.
$$
We denote by $[B]_{B_1}$ the $\sim_{B_1}$-equivalence class of a subalgebra $B$ in $\SA(B_1)$.

\begin{notation}
For $B$  in $\SA(B_1)$ let 
\begin{eqnarray*}
X(B_1,B_2;B) &=& \{ u \in \mathbb{U}(M_N): uB_2u^* \cap B_1 = B \}, \\
Y(B_2;B) &=&  \{ u \in \mathbb{U}(M_N) : u^*Bu \subseteq B_2  \},  \\
Z(B_1,B_2;[B]_{B_1}) &=& \{ u \in \mathbb{U}(M_N): uB_2u^*\cap B_1 \sim_{B_1} B \}. 
\end{eqnarray*}
\end{notation}

It is straightforward  that the complement of $\Delta(B_1,B_2)$ is  precisely the union of the sets $Z(B_1,B_2;[B]_{B_1})$, where  $B$ runs over all unital $C^*$-subalgebras of $B_1$ and $B\not= \mathbb{C}$ . Just for a moment, with out being  formal,  we may think $Z(B_1,B_2;[B]_{B_1})$ as being  parametrized by two coordinates. The first one is an algebra $\sim_{B_1}$-equivalent to $B$. Hence the first coordinate lives in $[B]_{B_1}$. The second, is a unitary $u$ that realizes the first coordinate as $uB_2u^*\cap B_1$.  $X(B_1,B_2;B)$ comes into play in order to parametrize this second coordinate. The problem is that $X(B_1,B_2;B_{B_1})$ is  complicated to handle (for instance it may not be closed). This is way we introduce the friendlier set $Y(B_2;B)$. Good properties about  $Y(B_2;B)$ is that it is a closed subset of $\mathbb{U}(M_N)$,  in fact  we will show it is a finite union of  enbedded compact submanifolds of $\mathbb{U}(M_N)$,  and it  contains  $X(B_1,B_2;B)$.

The rest of this  section  is the formalization of the previous idea. In concrete our first  goal  is to show $[B]_{B_1}$ has a structure of manifold  and we are particularly interested in finding  its dimension.

Let $\Sb(B_1,B)$ denote the $\sim_{B_1}$-stabilizer of $B$  i.e.
$$
\Sb(B_1,B)=\{ u \in \mathbb{U}(B_1) : uBu^*=B \}.
$$

\begin{remark}\label{StructureManifold[B]}
Given $B$ in $\SA(B_1)$ we can  endow $[B]_{B_1}$ with a structure of manifold. Indeed, let $\mathbb{U}(B_1)/ \Sb(B_1,B)$ denote the set of left-cosets and  consider the map
\begin{eqnarray*}
 \beta_{B}: [B]_{B_1} &\to &  \mathbb{U}(B_1)/ \Sb(B_1,B) , \\
 \beta_B (uBu^*)  &=&   u\Sb(B_1,B) .
\end{eqnarray*}
One can check $\beta_B$  is well defined and  bijective. Since $\mathbb{U}(B_1)/\Sb(B_1,B)$ is a manifold, $\beta_B$ induces a structure
of manifold on $[B]_{B_1}$. 
To avoid ambiguity we have to check the topology does not depend on the representative $B$. In fact, we will show the topology induced by $\beta_B$ is the same as the topology induced by the Hausdorff distance.

For $C_1$ and $C_2$ in $[B]_{B_1}$ define
$$
\dH(C_1,C_2)=\max \left\{   \sup_{x_2 } \inf_{x_1 } \{   \|  x_1-x_2  \|   \}  , \sup_{x_1 } \inf_{x_2 } \{   \|  x_1-x_2  \|   \}  \right\},
$$
where $x_i$ is taken in the unit ball of $C_i$, $i=1,2$.
Since unit balls of unital $C^*$-subalgebras of $B_1$ are compact subsets (in the norm topology), $\dH$ defines a metric on $[B]_{B_1}$.
Let $\tau$ and $\tau_H$ denote, respectively, the topologies on $[B]_{B_1}$ induced by $\beta_B$ and $\dH$. We are going to show $\tau=\tau_H$.  Consider the identity map $\id:([B]_{B_1}, \tau)\to ([B]_{B_1},\tau_H)$.
First we show $\id$ is continuous.
Since $\mathbb{U}(B_1)/\Sb(B_1,B)$ is endowed with the pull back topology from the quotient map
$\pi:\mathbb{U}(B_1)\to \mathbb{U}(B_1)/\Sb(B_1,B)$ where $\mathbb{U}(B_1)$ is taken with the norm topology,
$\id$ is continuous if and only if  the map 
$$
\beta_B^{-1} \circ\pi : \mathbb{U}(B_1) \to ([B]_{B_1},\tau_H)
$$
is continuous.
Take $(u_n)_{n\geq1}$ a sequence in $\mathbb{U}(B_1)$ and a unitary $u$ in $\mathbb{U}(B_1)$ such that $\lim_n \| u_n -u \| =0$. We need to show 
$$
\lim_n \dH( \beta_B^{-1} \circ\pi(u_n),\beta_B^{-1} \circ\pi(u)  )=\lim_n \dH( u_nBu_n^*,uBu^*  )=0.
$$
Take  $n_0$ such that $\|u_n -u \| < \varepsilon/2$ for all $n\geq n_0$.
For any $b$ in the unit ball of $B$ and any $n \geq n_0$, we have
$$
\|u_nbu_n^*-ubu^*\| < \varepsilon.
$$
Thus, for $n\geq n_0$ 
$$
\sup_{x_2}\inf_{x_1} \|x_1 -x_2\| < \varepsilon
$$
and 
$$
\sup_{x_1}\inf_{x_2} \|x_1 -x_2\| < \varepsilon,
$$
where $x_2$ is taken in the unit ball of $u_nBu_n^*$ and $x_1$ is taken in the unit ball of $uBu^*$.
Hence $\id:([B]_{B_1},\tau)\to ([B]_{B_1},\tau_H)$ is continuous. Lastly, since $\id$ is bijective, $([B]_{B_1},\tau)$ is compact and $([B]_{B_1},\tau_H)$ is Hausdorff, we conclude that $\id$ is a homeomorphism. Thus $\tau=\tau_H$. 
\end{remark}

Now that we know $[B]_{B_1}$ is a manifold, we want to find its dimension. Since by construction $[B]_{B_1}$ is diffeomorphic to $\mathbb{U}(B_1)/\Sb(B_1,B)$, $\dim [B]_{B_1}=\dim \mathbb{U}(B_1)- \dim \Sb(B_1,B)$. Thus we only need to find $\dim \Sb(B_1,B)$.

\begin{notation}
Whenever we take commutators they will be with respect to the ambient algebra $M_N$, in  other words for a subalgebra $A$ in $\SA(M_N)$
$$
A'=\{ x \in M_N: xa=ax, \quad \textrm{for all  $a$ in $A$} \}.
$$
Recall that $C(A)$ denotes the center of $A$ i.e.
$$
C(A)=A\cap A'=\{a\in A: xa=ax \quad \textrm{for all $x$ in A} \}.
$$
\end{notation}

\begin{prop}\label{DimOfStab}
For any $B_1$ in $\SA(M_N)$ and for any  $B$ in $\SA(B_1)$, we have
$$
\dim \Sb(B_1,B)= \dim \mathbb{U}(B) + \dim \mathbb{U}(B_1\cap B' ) - \dim  \mathbb{U}(C(B)).
$$
\end{prop}
\begin{proof}
We'll find a normal subgroup of $\Sb(B_1,B)$, for  which we can compute  its dimension and that partitions $\Sb(B_1,B)$  into a finite number of cosets.
Let $G$ denote the subgroup of $\Sb(B_1,B)$ generated by  $\mathbb{U}(B_1 \cap B')$ and $\mathbb{U}(B)$.  
Since the elements of $\mathbb{U}(B)$ commute with
the elements of  $\mathbb{U}(B_1\cap B')$, a typical element of $G$ looks like $vw$, where $v$ lies in  $\mathbb{U}(B)$ and $w$ lies  in $ \mathbb{U}(B_1\cap B')$. Taking into account  compactness of $\mathbb{U}(B)$ and $\mathbb{U}(B_1 \cap B')$, we deduced $G$ is compact.  

Now we show $G$ is normal in  $\Sb(B_1,B)$. Take $u$ an element  in $\Sb(B_1,B)$.  For a unitary  $v$ in $\mathbb{U}(B) $ it is immediate that $uvu^* $ lies  in $\mathbb{U}(B) $. For a unitary $w$ in $\mathbb{U}(B_1\cap B')$, the following computation shows $uwu^*$ belongs  to  $ \mathbb{U}(B_1 \cap  B ')$.
For any element $b$ in $B$ we have:
\begin{eqnarray*}
(uwu^* )b=uw(u^*bu)u^* =u(u^*bu)wu^*=b(uwu^*),
\end{eqnarray*}
where in the second equality we  used $u^*bu$ lies in  $ B$.
In conclusion $uGu^*$ is contained in  $G$ for all $u$ in $St(B_1,B)$ i.e. $G$ is normal in $\Sb(B_1,B)$.
 
As a result $\Sb(B_1,B)/G$ is a Lie group. The next  step  is to show  $\Sb(B_1,B)/G$ is finite.
Decompose $B$ as
$$
B= \oplus_{i=1}^I \oplus_{j=1}^{J_i}B(i,j) ,
$$ 
where for all $i$ there is $k_i$ such that  for $1\leq j \leq J_i$, $B(i,j)$ is $*$--isomorphic to $M_{k_i}$. For the rest  of our proof
we fix a  family, $ \beta(i,j): B(i,j) \to M_{k_i} $, of $*$--isomorphisms.

An element $u$ in $\Sb(B_1,B)$ defines a $*$--automorphism of $B$ by conjugation.
As a consequence, Propositions  \ref{2ndDecompositionAutomorphisms} and \ref{ThirdDecompositionAutomorphisms}  imply there are permutations $\sigma_i $ in $S_{J_i}$ and unitaries $v_i$ in $\mathbb{U}(\oplus_{j=1}^{J_i}B(i,j))$ such that 
\begin{equation}\label{TracialReduction}
\forall  b\in B : ubu^*= v \psi(b) v^*
\end{equation}
where $v=\oplus_{i=1}^I v_i $ is a uitary in $\mathbb{U}(B)$ and $\psi = \oplus_{i=1}^I \psi_{\sigma_i} $ is a $*$--automorphism in $Aut(B)$ (the maps $\psi$  depends on  the family of $*$--isomorphisms $\beta(i,j)$ we fixed earlier).
Equation (\ref{TracialReduction}) is telling  us important information. Firstly, that $\psi$ extends to an $*$--isomorphism of $B_1$ and most importantly,  this extension is an inner $*$--automorphism. Fix a unitary $U_{\psi}$ in $\mathbb{U}(B_1)$ such that $\psi(b)=AdU_{\psi}(b)$ for all $b$ in $B$ (note that $U_{\psi}$ may not be unique but we just pick one and fix it  for rest of the proof ). 
From equation (\ref{TracialReduction}) we deduce there is a unitary $w$ in $\mathbb{U}(B_1\cap B')$ satisfying  $u=vU_{\psi} w$. Since the number of functions $\psi$, that may arise from (\ref{TracialReduction}), is at most $J_1! \cdots J_I !$, we conclude
$$
|\Sb(B_1,B)/G| \leq J_1! \cdots J_I !
$$
Now that we know $\Sb(B_1,B)/G$ is finite we have $\dim \Sb(B_1,B)=\dim G$,  and Corollary \ref{ProductOfLieSubgroups} gives the result.
\end{proof}

From Proposition \ref{DimOfStab} and Remark \ref{StructureManifold[B]}, we get the following corollary.
\begin{coro}\label{Dim[B]}
For any $B_1$ in $\SA(M_N)$  and any  $B$ in $\SA(B_1)$, we have
$$
\dim [B]_{B_1}=\dim \mathbb{U}(B_1) - \dim \mathbb{U}(B'\cap B_1)
+ \dim \mathbb{U}(C(B)) - \dim \mathbb{U}(B)$$
\end{coro}

Now we focus our efforts on $Y(B_2;B)$.

\begin{prop}\label{DimYB}
Assume $Y(B_2;B) \not= \emptyset$. Then $Y(B_2;B)$ is a finite disjoint union of  embedded submanifolds of $\mathbb{U}(M_N)$. For each one of these submanifolds  there is $u \in Y(B_2;B)$ such that the submanifold's dimension is
\begin{eqnarray*}
\dim \Sb(M_N,B) + \dim \mathbb{U}(B_2)-\dim   \Sb(B_2,u^*Bu).
\end{eqnarray*}
Using Proposition \ref{DimOfStab} the later  equals
\begin{eqnarray}\label{DimYBeqn}
\dim \mathbb{U}(B') + \dim \mathbb{U}(B_2)-\dim \mathbb{U}(B_2\cap u^*B'u).
\end{eqnarray}

\end{prop}

\begin{proof}

We'll define an action on $Y(B_2;B)$ which will partition $Y(B_2;B)$ into a finite number of orbits,  each orbit an embedded submanifold of dimension  (\ref{DimYBeqn}) for a corresponding unitary.
Define an action of $\Sb(M_N,B) \times \mathbb{U}(B_2)$ on $Y(B_2;B)$ via
$$
(w,v).u=wuv^*.
$$ 
For $u \in Y(B_2;B)$  let  $\mathcal{O}(u)$ denote the orbit of $u$ and let $\mathcal{O}$ denote the set of all orbits.  
To prove $\mathcal{O}$ is finite  consider the function
\begin{eqnarray*}
\varphi: \mathcal{O}  &\to & \SA(B_2)/\sim_{B_2} , \\
\varphi(\mathcal{O}(u)) & = & [u^*Bu]_{B_2}.
\end{eqnarray*}
Firstly, we  need to show $\varphi$ is well defined. Assume $u_2 \in \mathcal{O}(u_1)$ and take $(w,v) \in \Sb(M_n,B) \times \mathbb{U}(B_2)$ such that $u_2=wu_1v^*$. From the identities
\begin{eqnarray*}
u_2^*Bu_2=vu_1w^*Bwu_1v^*=vu_1Bu_1v^*
\end{eqnarray*}
we obtain $[u_2Bu_2^*]_{B_2}=[u_1Bu_1^*]_{B_2}$. Hence $\varphi$ is well defined.

The next step is to show $\varphi$ is injective. Assume $\varphi(\mathcal{O}(u_1))=\varphi(\mathcal{O}(u_2))$, for $u_1,u_2 \in Y(B_2;B)$. Since $[u_1^*Bu_1]_{B_2}=[u_2^*Bu_2]_{B_2}$, we have $u_2^*Bu_2=vu_1^*Bu_1v^*$  for some  $v \in \mathbb{U}(B_2)$. But this implies $u_1v^*u_2^*\in \Sb(M_N,B)$ so if $w=u_1v^*u_2^*$ we conclude $(w,v).u_2=u_1$ which yields $\mathcal{O}(u_1) = \mathcal{O}(u_2)$. 
We conclude $|\mathcal{O}| \leq |\SA(B_2) / \sim_{B_2}| < \infty$.

Now we prove   each orbit is an embedded submanifold of $\mathbb{U}(M_N)$ of  dimension  (\ref{DimYBeqn}).
Since $\Sb(M_n,B)\times \mathbb{U}(B_2)$ is compact,  every orbit $\mathcal{O}(u)$ is compact. Thus, Proposition \ref{Orbit&Stabilizer}   implies $\mathcal{O}(u)$ is an embedded submanifold of $\mathbb{U}(M_N)$, diffeomorphic to
$$
( \Sb(M_N,B) \times \mathbb{U}(B_2) )/ \Sb(u)
$$
where
$$
\Sb(u)=\{  (w,v) \in \Sb(M_N,B) \times \mathbb{U}(B_2):  (w,v).u=u \}.
$$
Since
\[
(w,v).u=u\quad\Leftrightarrow\quad wuv^*=u\quad\Leftrightarrow\quad u^*wu=v,
\]
we deduce the group $\Sb(u)$ is  isomorphic to 
$$\mathbb{U}(B_2)\cap [  u^*\Sb(M_N,B)u ], $$
via the map $(w,v)\mapsto v$.
A straightforward computation shows
\[
u^*\Sb(M_N,B)u=\Sb(M_N,u^*Bu),
\]
for any $u\in \mathbb{U}(M_N)$. Hence, for any $u\in Y(B_2;B)$,
$$
\dim \mathcal{O}(u)=\dim \Sb(M_N,B) + \mathbb{U}(B_2) - \dim \mathbb{U}(B_2)\cap  \Sb(M_N,u^*Bu).
$$
Lastly, one can check $\mathbb{U}(B_2)\cap \Sb(M_N,u^*Bu)=\Sb(B_2,u^*Bu)$.
\end{proof}

\begin{notation}
For a unital $C^*$-subalgebra $B$ of  $B_1$, with the property that $B$ is unitarily equivalent to a $C^*$-subalgebra of $B_2$, or in other words  $Y(B_2;B)$ is nonempty, define
$$
d(B):=\dim [B]_{B_1}+ \max_i \{ \dim Y_i(B_2;B) \},
$$
where  $Y_1(B_2,B),\dots, Y_r(B_2;B)$ are disjoint submanifolds of $\mathbb{U}(M_N)$ whose union is $Y(B_2;B)$.
\end{notation}

As we mention at the beginning of this section, in order to prove Theorem \ref{DensityOfSmallIntersection}, we need to parametrize each $Z(B_1,B_2;[B]_{B_1})$ with a number of coordinates less than $N^2$. The number of coordinates will be given by $d(B)$. Thus  the next step is to show that, under the hypothesis  of Theorem \ref{DensityOfSmallIntersection}, we have $d(B)<N^2$  for $B \not= \mathbb{C}.$
We will later see that it suffices
to show $d(B)<N^2$  for $B\not= \mathbb{C}$  and $B$ abelian.

Before we proceed, we recall definition of multiplicity of  of a representation. The following lemma combines   Lemma III.2.1 in ~\cite{Davidson} and Theorem 11.9 in ~\cite{Takesaki}.

\begin{lemma}\label{PartialMultiplicities}
Suppose $\varphi: A_1 \to A_2$ is a unital $*$-homomorphism and $A_i$ is isomorphic to $\bigoplus_{j=1}^{l_i}M_{k_i(j)}$, ($i=1,2$). Then $\varphi$ is determined, up to unitary equivalence in $A_2$, by an $l_2\times l_1$ matrix,  written
$\mu=\mu(\phi)=\mu(A_2,A_1)$, having nonnegative integer entries such that
$$
\mu\left[ \begin{array}{c}
k_1(1) \\
\vdots \\
k_1(l_1)
\end{array}
\right]=
\left[ \begin{array}{c}
k_2(1) \\
\vdots \\
k_2(l_2)
\end{array}
\right].
$$  
We call this the matrix of partial multiplicities. In the special case when $\varphi$ is a  unital $*$--representation  of $A_1$  into $M_N$, $\mu$ is a row vector  and this vector is called the  multiplicity of the representation.
One constructs $\mu$ as follows: decompose $A_p$ as
$$
A_p=\oplus_{j=1}^{l_p}A_p(j)
$$
where each $A_p(j)$ is  simple, $p=1,2$, $1\leq j \leq l_p$. Taking  projections, $\pi$ induces  unital $*$--representations $\pi_i:A_1 \to A_2(i)$, $1\leq i \leq l_2$. But up to unitary equivalence, $\pi_i$ equals
$$
\underbrace{\id_{A_1(1)} \oplus \cdots \oplus \id_{A_1(1)}}_{m_{i,1}-\textrm{times}}\oplus \cdots \oplus \underbrace{\id_{A_1(l_1)} \oplus \cdots \oplus \id_{A_1(l_1)}}_{m_{i,l_1}-\textrm{times}}
$$
for some nonnegative integer $m_{i,j}$, $1\leq j \leq l_1$. Set $\mu[i,j]:=m_{i,j}$.
In particular, $\mu[i,j]$ equals the rank of $\pi_i(p)\in A_2(i)$, where $p$ is a minimal projection in $A_1(j)$.
Clearly, $\pi$ is injective if and only if for all $j$ there is $i$ such that $\mu[i,j]\ne0$.

Furthermore,  the $C^*$-subalgebra
$$
A_2\cap\varphi(A_1)'=\{ x\in A_2: x\varphi(a)=\varphi(a)x \quad \textrm{for all $a\in A_1$} \}
$$
is $*$--isomorphic to
$
\bigoplus_{i=1}^{l_2}\bigoplus_{j=1}^{l_1} M_{\mu[i,j]}
$
and if we have morphisms $A_1\to A_2\to A_3$, then $\mu(A_3,A_2)\mu(A_2,A_1)=\mu(A_3,A_1)$ for the corresponding
matrices.
\end{lemma}

Our next task is to show $d(B)<N^2$, for abelian $B\not=C$. We prove it by cases, so let us start.

\begin{lemma}\label{EmbeddngIntoGCD}
Assume $B_i$ is  $*$--isomorphic to $M_{k_i}$,  ($i=1,2$) and let  $k= \gcd (k_1,k_2)$.
Take  $B $ a unital $C^*$-subalgebra of $B_1$  such that it is unitarily equivalent to a $C^*$-subalgebra of $B_2$. Then there is an injective unital $*$--representation of $B$ into $M_k$.
\end{lemma}
\begin{proof}
Take $u$ in $Y(B_2;B)$ so that $u^*Bu\subseteq B_2$.
Let $m_i:=\mu(M_N,B_i)$, so that $m_ik_i=N$, ($i=1,2$).
 Find positive integers $p_1$ and $p_2$ such that  $k_1=kp_1$ and $k_2=kp_2$
 Assume $B$  is $*$--isomorphic to 
$\bigoplus_{j=1}^l M_{n_j}$.
To prove the result it is enough  to show there are positive integers $(m(1),\dots m(l))$ such that
$$
n_1 m(1)+ \dots + n_l m(l) =k.
$$
Let 
\begin{eqnarray*}
\mu(B_1,B)=[m_1(1),\dots, m_1(l)],\\
\mu(B_2,u^*Bu)=[m_2(1),\dots ,m_2(l)].
\end{eqnarray*}
Since $\mu(M_N,B_1)\mu(B_1,B)=\mu(M_N,B_2)\mu(B_2,u^*Bu)$  we deduce that   
$m_1m_1(j)=m_2m_2(j)$ for all $1\leq j \leq l$. Multiplying by $k$ and using $N=m_1k_1=m_2k_2$ we conclude
$$
\frac{N}{p_1}m_1(j)=km_1m_1(j)=km_2m_2(j)=\frac{N}{p_2}m_2(j),
$$
so $p_2m_1(j)=p_1m_2(j)$.
Since  $\gcd(p_1,p_2)=1$,  the number $\tfrac{m_1(j)}{p_1}=\tfrac{m_2(j)}{p_2} $ is a positive  integer
whose value we name $m(j)$.
From
$$
kp_1=k_1 = \sum_{j=1}^l n_j m_1(j)=\sum_{j=1}^l n_j m(j)p_1,
$$
we conclude $k=\sum_{j=1}^l n_jm(j)$.
\end{proof}

\begin{prop}\label{dBLessN2SimpleCase}
Assume $B_1$ and $B_2$ are simple. Take $B\not= \mathbb{C}$ an abelian unital $C^*$-subalgebra of $B_1$,  that is unitarily equivalent to a $C^*$-subalgebra of $B_2$. Then $d(B) < N^2$.
\end{prop}
\begin{proof}
Assume $B_i$ is $*$--isomorphic to $M_{k_i}$, ($i=1,2$) and $B$ is $*$--isomorphic to $\mathbb{C}^l$, $l\geq 2$.
Using Corollary~\ref{Dim[B]} and Proposition~\ref{DimYB}, we may take
$u$ in $Y(B_2,B)$ such that  $d(B)$  equals the sum of the following terms,
\begin{eqnarray*}
S_1(B)&:= & \dim \mathbb{U}(B_1 )- \dim \mathbb{U}(B_1\cap B') , \\
 S_2(B)&:= &  \dim  \mathbb{U}(B_2 )- \dim \mathbb{U}(B_2\cap u^*B'u) , \\
 S_3(B)&:= & \dim \mathbb{U}(B') ,
\end{eqnarray*}
Let $k=\gcd(k_1,k_2)$ and write $k_1=kp_1$, $k_2=kp_2$.
From proof of Lemma \ref{EmbeddngIntoGCD}, there are positive integers $m(j)$, $1\leq j \leq l$, such that 
\begin{align*}
\mu(B_1,B)&=[m(1)p_1,  \dots , m(l)p_1] \\
\mu(B_2,B)&=[m(1)p_2,  \dots , m(l)p_2].
\end{align*}
Hence 
\begin{align*}
S_1(B)&=k_1^2 - \sum_{i=1}^l  m(i)^2p_1^2= k^2p_1^2 - \sum_{i=1}^l m(i)^2 p_1^2 \\
S_2(B)&=k_2^2 - \sum_{i=1}^l m(i)^2p_2^2= k^2p_2^2 - \sum_{i=1}^l m(i)^2 p_2^2 .
\end{align*}
Let $m_i=\mu(M_N,B_i)$, ($i=1,2$).
Since
\[
\mu(M_N,B_1)\mu(B_1,B)=\mu(M_N,B_2)\mu(B_2,u^*Bu),
\]
 we get
\begin{equation}\label{EqualityMulRepFromB}
\mu(M_N,B)
\begin{aligned}[t]
&=[m_1p_1m(1),  \dots , m_1p_1m(l) ] \\
&=[m_2p_2m(1), \dots , m_2p_2m(l)].
\end{aligned}
\end{equation}
Hence 
$$
S_3(B)=\sum_{i=1}^l (m(i)p_1m_1)(m(i)p_2m_2)=\left( \sum_{i=1}^l m(i)^2\right)p_1p_2m_1m_2.
$$
Factoring the term $\sum_{i=1}^l m(i)^2$  we get $d(B)$ equals
$$
\left( \sum_{i=1}^l m(i)^2 \right)\left(  p_1p_2m_1m_2-p_1^2-p_2^2  \right)+ k^2(p_1^2 + p_2^2).
$$
On the other hand, using $N=m_1k_1=m_1kp_1=m_2k_2=m_2kp_2$, we get $N^2=k^2p_1p_2m_1m_2$. 
Hence $d(B)< N^2$ if and only if
\begin{equation}\label{ReductiondBlessN2}
\left( \sum_{i=1}^l m(i)^2 \right)\left(  p_1p_2m_1m_2-p_1^2-p_2^2  \right)<  k^2(p_1p_2m_1m_2 -p_1^2 - p_2^2).
\end{equation}
We want to cancel $(p_1p_2m_{1}m_{2} -p_1^2 - p_2^2), $  in   equation  (\ref{ReductiondBlessN2}), so we  prove it is positive. First we divide it by $p_1p_2$ to get  $m_1m_2-\frac{p_1}{p_2} - \frac{p_2}{p_1}$. But from equation $(\ref{EqualityMulRepFromB})$ we have $\frac{p_1}{p_2}=\frac{m_2}{m_1}$. Thus  we need to show $m_1m_2-\frac{m_1}{m_2} - \frac{m_2}{m_1}$ is positive. If we divide it by $m_{1}m_{2}$ we get $1-\frac{1}{m_1^2}-\frac{1}{m_2^2}$, which is clearly positive
(recall that  $m_1 \geq 2$ and $m_2 \geq 2$ since $B_1 \not= M_N$ and $B_2 \not= M_N$).
Therefore, equation (\ref{ReductiondBlessN2}) is equivalent to
$$
\sum_{i=1}^l m(i)^2 < k^2.
$$
But $\sum_{i=1}^l m(i)=k$, $l\geq 2$ and each $m(i)$ is positive.
\end{proof}

 In the nonsimple case in Theorem \ref{DensityOfSmallIntersection}, we will need  some minimization lemmas  to show $d(B)<N^2$, for abelian $B\not=\mathbb{C}$. 
A straightfroward use of Lagrange multipliers proves the following lemma, and the one after that is even more elementary.
\begin{lemma}\label{MinimizationLagrange}
Fix a positive integer $n$ and let $r_1,\dots , r_n$ be positive real numbers. Then
\[ 
\min\bigg\{  \sum_{j=1}^n \frac{x_j^2}{r_j} \;\bigg|\; \sum_{j=1}^n x_j=1  \bigg\} =\frac{1}{\sum_{j=1}^n r_j},
\] 
where the minimum is taken over all $n$-tuples of real numbers that sum up to 1.
\end{lemma}

\begin{lemma}\label{MaximizationLemma}
For  an integer $k\geq 2$ define
$$
h(x,y)=2xy-\bigg( 1+\frac{1}{k^2}\bigg)y^2-\frac{1}{2}x^2.
$$
Then
$$
\max\{  h(x,y) \mid 0\leq x \leq 1,\, 0 \leq y \leq 1/2\}  =\frac{1}{4}-\frac{1}{4k^2}.
$$
\end{lemma}

\begin{prop}\label{dBLessN2NonSimpleCase}
Suppose $\dim C(B_1) \geq 2$  and  $B_1$ is $*$--isomorphic to
\begin{equation}\label{eq:B1}
M_{N/\dim C(B_1)}\oplus \cdots \oplus M_{N/\dim C(B_1)}.
\end{equation}
Assume one of the following cases holds:
\begin{enumerate}
\item $\dim C(B_2)=1$,
\item $B_1$ is $*$--isomorphic to
$$
M_{N/2} \oplus M_{N/2},
$$
$B_2$ is $*$--isomorphic to 
$$
M_{N/2}\oplus M_{N/(2k)}
$$
where $k \geq 2$.
\item $\dim C(B_2)\geq 3$  and $B_2$ is $*$--isomorphic to
$$
M_{N/\dim C(B_2)}\oplus \cdots \oplus M_{N/\dim C(B_2)}.
$$
\end{enumerate}
Then for any $B\not=\mathbb{C}$ an abelian unital $C^*$-subalgebra of $B_1$ that is unitarily equivalent to a $C^*$-subalgebra of $B_2$, we have that $d(B)<N^2$.
\end{prop}
\begin{proof}
Let $l_i=\dim C(B_i)$, ($i=1,2$), $l=\dim(B)$.
Take $u$ in $Y(B_2;B)$ such that $d(B)$ is the sum of the following terms:
\begin{eqnarray}
S_1(B)&:=& \dim \mathbb{U}(B_1 )- \dim \mathbb{U}(B_1\cap B') ,  \\
 S_2(B)&:= &  \dim  \mathbb{U}(B_2 )- \dim \mathbb{U}(B_2\cap u^*B'u) , \\
 S_3(B)&:= & \dim \mathbb{U}(B') .
 \end{eqnarray}
Write
\begin{eqnarray*}
\mu(B_1,B)&=&[a_{i,j}]_{1\leq i \leq l_1, 1\leq j \leq l}, \\
\mu(B_2,u^*Bu)&=&[b_{i,j}]_{1\leq i \leq l_2, 1\leq j \leq l},\\
\mu(M_N,B_1)&=&[m_1(1), \dots , m_1(l_1)],\\
\mu(M_N,B_2)&=&[m_2(1), \dots , m_2(l_2)],\\
\mu(M_N,B)&=&[m(1), \dots ,m(l)].
\end{eqnarray*}
Then
\begin{eqnarray*}
S_1(B)&=&\frac{N^2}{l_1}-\sum_{i=1}^{l}\sum_{j=1}^{l_1} a_{i,j}^2,\\
S_2(B)&=&\dim \mathbb{U}(B_2)-\sum_{i=1}^{l}\sum_{j=1}^{l_2} b_{i,j}^2,\\
S_3(B)&=&\sum_{j=1}^lm(j)^2.
\end{eqnarray*}
Since the sum of the ranks appearing in~\eqref{eq:B1} is $N$, we have
$m_1(i)=1$ for   all $1\leq i \leq l_1$.
Since \[
\mu(M_N,B)=\mu(M_N,B_1)\mu(B_1,B)=\mu(M_N,B_2)\mu(B_2,u^*Bu),
\]
we must have
$$
m(j)=\sum_{i=1}^{l_1}a_{i,j}=\sum_{i=1}^{l_2}m_2(i)b_{i,j}
$$
for all $1\leq j \leq l$.
Hence there are nonnegative numbers $\alpha_{i,j}$ and $\beta_{i,j}$ such that $\sum_{i=1}^{l_1}\alpha_{i,j}=\sum_{i=1}^{l_2}\beta_{i,j}=1$ and $a_{i,j}=\alpha_{i,j}m(j)$, $m_2(i)b_{i,j}=\beta_{i,j}m(j)$.
On the other hand, since $B$ is a unital $C^*$-subalgebra of $M_N$  we must have
$$\sum_{j=1}^lm(j)=N.$$
Thus, there are positive numbers $\gamma_j$, ($1\leq j \leq l$), such that $\sum_{j=1}^l \gamma_j =1$ and $m(j)=\gamma_jN$.
It will be important to notice that $\gamma_j >0$ for all $1\leq j \leq l$ ( otherwise $B$ is not a  unital $C^*$-algebra of $M_N$).
In consequence,
 \begin{eqnarray*}
S_1(B)&=&\frac{N^2}{l_1}-N^2\bigg( \sum_{j=1}^l  \gamma_j^2 \bigg( \sum_{i=1}^{l_1}\alpha_{i,j}^2 \bigg) \bigg),\\
S_2(B)&=&\dim\mathbb{U}(B_2)-N^2\bigg( \sum_{j=1}^l \gamma_j^2  \bigg(\sum_{i=1}^{l_2}\frac{\beta_{i,j}^2}{m_2(i)^2}\bigg) \bigg),\\
S_3(B)&=&N^2\bigg(\sum_{j=1}^l \gamma_j^2 \bigg).
\end{eqnarray*}

\noindent
{\em Case (1)}.
$B_2$ is simple, let us say  it is $*$--isomorphic to $M_{k_2}$. In this case $\mu(M_N,B_2)=[m_2]$ is just one number and we must have $m_2k_2=N$. Notice that $m_2 \geq2$, since by our standing assumption, $B_2\neq M_N$.
Also notice that from $\mu(M_N,B_2)\mu(B_2,u^*Bu)=\mu(M_N,B)$ we  obtain $m_2b_{i,1}=m(i)$ and $\beta_{i,1}=1$ for all $1\leq i \leq l$.
In consequence
 \begin{eqnarray*}
S_1(B)&=&\frac{N^2}{l_1}-N^2\bigg( \sum_{j=1}^l  \gamma_j^2 \bigg( \sum_{i=1}^{l_1}\alpha_{i,j}^2 \bigg) \bigg),\\
S_2(B)&=&\frac{N^2}{m_2^2}-\frac{N^2}{m_2^2}\bigg( \sum_{j=1}^l \gamma_j^2  \bigg),\\
S_3(B)&=&N^2\bigg(\sum_{j=1}^l \gamma_j^2 \bigg).
\end{eqnarray*}
From Lemma \ref{MinimizationLagrange}, we deduce
 \begin{eqnarray*}
S_1(B)& \leq &\frac{N^2}{l_1}-\frac{N^2}{l_1}\bigg( \sum_{j=1}^l  \gamma_j^2  \bigg).
\end{eqnarray*}
 Thus, it suffices to show
 $$
N^2\bigg(  \frac{1}{l_1}+\frac{1}{m_2^2} + \sum_{j=1}^{l}\gamma_j^2 -  \frac{1}{l_1}\bigg(\sum_{j=1}^{l}\gamma_j^2\bigg) -\frac{1}{m_2^2}\bigg( \sum_{j=1}^l \gamma_j ^2 \bigg) \bigg)< N^2
 $$ 
or equivalently 
$$
\bigg(\sum_{j=1}^{l}\gamma_j^2 \bigg) \bigg( 1 -\frac{1}{l_1}-\frac{1}{m_2^2} \bigg)< 1-\frac{1}{l_1}-\frac{1}{m_2^2} .
 $$
Since $l_1 \geq 2$ and $m_2 \geq 2$ we can cancel the term $1-\frac{1}{l_1}-\frac{1}{m_2^2}$. Thus we need to show $\sum_{j=1}^l \gamma_j^2<1$. But the latter  follows from the fact that $l \geq 2$, each $\gamma_j $ is positive and $\sum_{j=1}^l\gamma_j=1$.

\medskip
\noindent
{\em Case (2)}.
We have
\begin{eqnarray*}
\mu(M_N,B_1)=[1,1],\\
\mu(M_N,B_2)=[1,k].
\end{eqnarray*}

Thus
 \begin{eqnarray*}
S_1(B)&=&\frac{N^2}{2}-N^2\bigg( \sum_{j=1}^l  \gamma_j^2 \bigg( \alpha_{1,j}^2 + \alpha_{2,j}^2\bigg) \bigg),\\
S_2(B)&=&\frac{N^2}{4}+ \frac{N^2}{4k^2}-N^2\bigg( \sum_{j=1}^l \gamma_j^2  \bigg(\beta_{1,j}^2 + \frac{\beta_{2,j}^2}{k^2}\bigg) \bigg),\\
S_3(B)&=&N^2\bigg(\sum_{j=1}^l \gamma_j^2 \bigg).
\end{eqnarray*}
From Lemma \ref{MinimizationLagrange} we  obtain 
\begin{eqnarray*}
S_1(B) & \leq & \frac{N^2}{2}-\frac{N^2}{2}\bigg( \sum_{j=1}^l \gamma_j^2 \bigg).\\
\end{eqnarray*}
Thus, it suffices to show
$$
  \frac{1}{2}+\frac{1}{4}+\frac{1}{4k^2} + \sum_{j=1}^{l}\gamma_j^2 \bigg(\frac{1}{2} -\beta_{1,j}^2-\frac{1}{k^2}\beta_{2,j}^2 \bigg) < 1
$$
or, equivalently,
\[ 
\sum_{j=1}^{l}\gamma_j^2 \bigg(\frac{1}{2} -\beta_{1,j}^2-\frac{1}{k^2}\beta_{2,j}^2 \bigg) < \frac{1}{4}-\frac{1}{4k^2}.
\] 
Define
\begin{equation}\label{eq:ro}
r= \sum_{j=1}^{l}\gamma_j^2 \bigg(\frac{1}{2} -\beta_{1,j}^2-\frac{1}{k^2}\beta_{2,j}^2 \bigg) .
\end{equation}
Now we use the constraints on the variables $\gamma_j$ and $\beta_{i,j}$.
First of all we have $\beta_{1,j}+\beta_{2,j}=1$ for all $1\leq i \leq l$.
Thus, $r $ simplifies to
$$
r=\sum_{j=1}^l \gamma_j^2 \bigg( 2\beta_{2,j}-\bigg( 1+\frac{1}{k^2} \bigg)\beta_{2,j}^2 -\frac{1}{2} \bigg).
$$
We also have
\begin{eqnarray}\label{ConstrainsOverGammaBeta}
\sum_{j=1}^l \beta_{2,j}\gamma_j&=&\frac{1}{2}.
\end{eqnarray}
Indeed, since all blocks of $B$ are one dimensional,  we must have
$$
\sum_{j=1}^l b_{2,j}=\frac{N}{2k}.
$$
But $kb_{2,j}=\beta_{2,j}m(j)=\beta_{2,j}\gamma_{j}N$, which implies (\ref{ConstrainsOverGammaBeta}).
The final constraint is $\sum_{j=1}\gamma_j=1$.

Now we make the  change of variables $q_j:=\gamma_j\beta_{2,j}$ and $r$ becomes
$$
r=2\bigg( \sum_{j=1}^l q_j \gamma_j \bigg) -\bigg( 1+\frac{1}{k^2}\bigg)\bigg(\sum_{j=1}^l q_j^2 \bigg)-\frac{1}{2}\bigg( \sum_{j=1}^l \gamma_j^2 \bigg).
$$
Letting $\gamma=(\gamma_1,\dots , \gamma_l)$ and $q=(q_1,\dots, q_l)$ and using the Cauchy-Schwartz inequality,
we get
$$
r\leq 2\| q\|_2\|\gamma\|_2 -\bigg( 1+\frac{1}{k^2} \bigg)\|q\|_2^2-\frac{1}{2}\|\gamma\|_{2}^2
$$
Set $x=\| \gamma\|$, $y=\|q\|$. Notice that $0\leq x \leq 1$ and $0\leq y \leq 1/2$. Take
$$h(x,y)=2xy-\bigg( 1+\frac{1}{k^2}\bigg)y^2-\frac{1}{2}y^2$$
apply Lemma \ref{MaximizationLemma} to get
\[ 
r\leq h(\|\gamma\|,\|q\|)\leq \frac{1}{4}-\frac{1}{4k^2}.
\] 
Now we will rule out equality. Assuming, for contradiction, $r=\frac14-\frac1{4k^2}$,
we must have equality in the instince of the Cauchy-Schwartz inequality. Hence $q=z\gamma$ for some real number $z$. Summing over the coordinates we deduce $z=1/2$ and  then, for all $1\leq j \leq l$,
$$\frac{1}{2}\gamma_j = q_j=\gamma_j \beta_{2,j} .$$
Since $\gamma_j>0$ we can cancel and get $\beta_{2,j}=1/2$.  Thus, using the original formulation~\eqref{eq:ro} of $r$, we get
$$
r=\bigg( \frac{1}{4}-\frac{1}{4k^2}\bigg)\bigg( \sum_{j=1}^l \gamma_j^2 \bigg) 
$$
which is strictly less that $1/4-1/(4k^2)$,  because $k\geq 2$, $l\geq2$, all $\gamma_j$ are strictly  positive and $\sum_{j=1}^l\gamma_j=1$. 

\medskip
\noindent
{\em Case (3)}.
Then $B_2$ is $*$--isomorphic to
$$
\underbrace{M_{N/l_2}\oplus \cdots \oplus M_{N/l_2}}_{l_2-\textrm{times}}.
$$
Arguing as we did before for $m_1(i)$, we have  $m_2(i)=1$, for all   $1\leq i \leq l_2$. 
Hence 
 \begin{eqnarray*}
S_1(B)&=&\frac{N^2}{l_1}-N^2\bigg( \sum_{j=1}^l  \gamma_j^2 \bigg( \sum_{i=1}^{l_1}\alpha_{i,j}^2 \bigg) \bigg),\\
S_2(B)&=&\frac{N^2}{l_2}-N^2\bigg( \sum_{j=1}^l \gamma_j^2  \bigg(\sum_{i=1}^{l_2}\beta_{i,j}^2\bigg) \bigg),\\
S_3(B)&=&N^2\bigg(\sum_{j=1}^l \gamma_j^2 \bigg).
\end{eqnarray*}
 From Lemma \ref{MinimizationLagrange} we deduce
 \begin{eqnarray*}
S_1(B)& \leq &\frac{N^2}{l_1}-\frac{N^2}{l_1}\bigg( \sum_{j=1}^l  \gamma_j^2  \bigg),\\
S_2(B)&\leq &\frac{N^2}{l_2}-\frac{N^2}{l_2}\bigg( \sum_{j=1}^l \gamma_j ^2  \bigg).
\end{eqnarray*}
 Thus, it suffices to show
 $$
N^2\bigg(  \frac{1}{l_1}+\frac{1}{l_2} + \sum_{j=1}^{l}\gamma_j^2 -  \frac{1}{l_1}\bigg(\sum_{j=1}^{l}\gamma_j^2\bigg) -\frac{1}{l_2}\bigg( \sum_{j=1}^l \gamma_j ^2 \bigg) \bigg)< N^2
 $$ 
or equivalently 
$$
\bigg(\sum_{j=1}^{l}\gamma_j^2 \bigg) \bigg( 1 -\frac{1}{l_1}-\frac{1}{l_2} \bigg)< 1-\frac{1}{l_1}-\frac{1}{l_2} .
 $$
Since $l_1 \geq 2$ and $l_2 \geq 3$ we can cancel the term $1-\frac{1}{l_1}-\frac{1}{l_2}$ in the above equation and finish the proof as in the previous case. 
 \end{proof}

The next step is to find parameterizations of  $Z(B_1,B_2;[B]_{B_1})$.

\begin{lemma}\label{BSimpleCIsomC2}
Take $B\not= \mathbb{C}$ a unital $C^*$-subalgebra of  $B_1$ that is unitarily equivalent to a $C^*$-subalgebra of $B_2$.
If $\dim \mathbb{U}(B_1)+\dim \mathbb{U}(B_2) \leq N^2$, $B$ is simple and  $C $ in $\SA(B)$ is $*$--isomorphic to $\mathbb{C}^2$, then $d(B) \leq d(C)$.
\end{lemma}
\begin{proof}
Assume $B$ is $*$--isomorphic to $M_k$ and let $m$ denote the multiplicity of $B$ in $M_N$. Thus we must have $km=N$.
Take a unitary $u$ in the submanifold of maximum dimension in $Y(B_2;B)$, so that $d(B)$ is the sum of the terms
\begin{eqnarray*}
S_1(B)&:=& \dim \mathbb{U}(B_1 )- \dim \mathbb{U}(B_1\cap B') , \\
 S_2(B)&:= &  \dim  \mathbb{U}(B_2 )- \dim \mathbb{U}(B_2\cap u^*B'u) ,\\
 S_3(B)&:= & \dim \mathbb{U}(B') , \\
 S_4(B)&:= & \dim \mathbb{U}(B\cap B')-\dim \mathbb{U}(B) .
\end{eqnarray*}
and let $v$ lie in the submanifold of maximum dimension in $Y(B_2,C)$ so that $d(C)$ is the sum of the terms
\begin{eqnarray*}
S_1(C)&:=& \dim \mathbb{U}(B_1 )- \dim \mathbb{U}(B_1\cap C') , \\
 S_2(C)&:= &  \dim  \mathbb{U}(B_2 )- \dim \mathbb{U}(B_2\cap v^*C'v) ,\\
 S_3(C)&:= & \dim \mathbb{U}(C') . \\
\end{eqnarray*}
Clearly, $S_4(B)=1-k^2$.
We write
\begin{eqnarray*}
B_1 & \simeq  & \bigoplus_{i=1}^{l_1}M_{k_1(i)},\\
B_2 & \simeq  & \bigoplus_{i=1}^{l_2}M_{k_2(i)}.
\end{eqnarray*}
and
\begin{eqnarray*}
\delta(B_1)&=&[k_1(1),\dots , k_1(l_1)]^t, \\
\delta(B_2)&=&[k_2(1),\dots , k_2(l_2)]^t.
\end{eqnarray*}
 From definition of multiplicity and the fact that it is invariant under unitary equivalence we get
\begin{eqnarray}
\mu(B_1,B)k&=&\delta(B_1),  \label{M(B1,B)DB} \\
 \mu(B_2,u^*Bu)k&=&\delta(B_2), \notag \\ 
\mu(M_N,B_1)\delta(B_1)&=&\mu(M_N,B_2)\delta(B_2)=N, \notag \\ 
 \mu(M_N,B_1)\mu(B_1,B) &=& \mu(M_N,B_2)\mu(B_2,u^*Bu)=m . \notag 
\end{eqnarray}
From Lemma \ref{PartialMultiplicities} and  equation (\ref{M(B1,B)DB}) we get
\begin{equation}\label{NiceFaceCommB}
\dim \mathbb{U}(B_1\cap B')=\frac{1}{k^2}\dim \mathbb{U}(B_1).
\end{equation}
Hence
$$
S_1(B)=\bigg( 1-\frac{1}{k^2} \bigg)\dim \mathbb{U}(B_1).
$$
Similarly
$$
S_2(B)=\bigg( 1-\frac{1}{k^2} \bigg)\dim \mathbb{U}(B_2).
$$

Now it is the turn of $C$. To ease notation let 
$$\mu(B,C)=[x_1, x_2]$$
Notice that $x_1+x_2=k$.
We claim 
$$
S_1(C)=\bigg(1- \frac{x_1^2+x_2^2}{k^2} \bigg) \dim \mathbb{U}(B_1).
$$
Using  $\mu(B_1,C)=\mu(B_1,B)\mu(B,C)$ we get
$$
\dim \mathbb{U}(B_1\cap C')= (x_1^2+x_2^2)\dim \mathbb{U}(B_1\cap B').
$$
Furthermore using (\ref{NiceFaceCommB}) we obtain
$$
\dim \mathbb{U}(B_1\cap C')= \frac{x_1^2+x_2^2}{k^2}\dim \mathbb{U}(B_1).
$$
Hence our claim follows from definition of $S_1(C)$.
Similarly
$$
S_2(C)=\bigg(1- \frac{x_1^2+x_2^2}{k^2} \bigg) \dim \mathbb{U}(B_2).
$$
Lastly from   $\mu(M_N,C)=[m x_1,m x_2]$ and $mk=N$  we get
\begin{eqnarray*}
S_3(C)&=&(x_1^2+x_2^2)\frac{N^2}{k^2} ,\\
S_3(B)&=&\frac{N^2}{k^2}.
\end{eqnarray*}
To prove $d(B) \leq d(C)$ we'll  show 
\begin{equation}\label{Reaccomodation}
S_1(B)-S_1(C)+S_2(B)-S_2(C)+S_4(B) \leq S_3(C)-S_3(B) .
\end{equation}
Using the description of each summand we have that left hand side of (\ref{Reaccomodation}) equals
\begin{eqnarray*}
 \frac{x_1^2+x_2^2-1}{k^2}\bigg(\dim \mathbb{U}(B_1)+ \dim \mathbb{U}(B_2) \bigg)  + 1-k^2.
\end{eqnarray*}
The  right hand side of (\ref{Reaccomodation}) equals
$$
\frac{x_1^2+x_2^2-1}{k^2}N^2.
$$
But $x_1$ and $x_2$ are strictly positive, because $C$ is a unital subalgebra of $B$.
Hence we can cancel $x_1^2+x_2^2-1$ and  finish  the proof by using that $1-\delta(B)^2<0$ and the assumption $\dim \mathbb{U}(B_1)+\dim \mathbb{U}(B_2) \leq  N^2$.
\end{proof}

We recall an important perturbation result that can  be found in  Lemma  III.3.2 from  ~\cite{Davidson}.
\begin{lemma}\label{PerturbationDavidson}
Let $A$ be a finite dimensional $C^*$-algebra. Given any positive number $\varepsilon$ there is a positive  number $\delta=\delta(\varepsilon)$ so that whenever $B$ and $C$ are unital $C^*$-subalgebras of  $A$ and such that $C$ has a system of matrix units $\{ e_C(s,i,j) \}_{s,i,j}$, satisfying $\dist(e_C(s,i,j),B)< \delta$ for all $s,i$ and $j$, then there is a unitary $u$ in $\mathbb{U}(C^*(B,C))$ with $\| u-1 \| < \varepsilon$ so that $uCu^* \subseteq B$.
\end{lemma}

\begin{notation}
For an element $x$ in  $M_N$ and a positive number $\varepsilon$, $\mathcal{N}_{\varepsilon}(x)$ denotes the open $\varepsilon$-neighborhood around $x$ (i.e. open ball of radius $\varepsilon$ centered at $x$), where the distance is  from the operator norm in $M_N$.
\end{notation}

The next proposition is quite technical and is mainly a  consequence of  Lemma \ref{PerturbationDavidson}.
The set $[B]_{B_1}$ is endowed with the equivalent topologies described in Remark~\ref{StructureManifold[B]}.
\begin{lemma}\label{ContinuityIntersectionProp}
Take $B$ in $\SA(B_1)$ and assume $Z(B_1,B_2;[B]_{B_1})$ is nonempty. 
Then the function 
\begin{eqnarray}\label{ContinuityIntersectionFunct}
Z(B_1,B_2;[B]_{B_1}) &  \to &   [B]_{B_1} \\
u & \mapsto &  uB_2u^* \cap B_1 \nonumber
\end{eqnarray}
is  continuous.
\end{lemma}
\begin{proof}
Assume $B$ is $*$--isomorphic to
$$
\bigoplus_{s=1}^l M_{k_s}.
$$
First we recall that the topology of $[B]_{B_1}$ is induced by the bijection 
\begin{eqnarray*}
\beta:[B]_{B_1} & \to &  \mathbb{U}(B_1)/\Sb(B_1,B), \\
\beta(uBu^*) & = & u\Sb(B_1,B).
\end{eqnarray*}
For convenience let $\pi:\mathbb{U}(B_1) \to \mathbb{U}(B_1)/\Sb(B_1,B)$  denote the canonical quotient map. 
Pick $u_0$ in $Z(B_1,B_2;[B]_{B_1})$. With no loss of generality  we may assume $B=u_0B_2u_0^* \cap B_1$. 

We prove the result by contradiction. 
Suppose the function in  (\ref{ContinuityIntersectionFunct}) is not continuous at $u_0$.
Then there is a sequence $(u_k)_{k \geq 1} \subset Z(B_1,B_2,[B]_{B_1})$ 
and an open neighborhood  $\mathcal{N}$ of $B$ in $[B]_{B_1}$ 
such that
\begin{enumerate}
\item $\lim_k  u_k = u_0$,
\item for all $k$,  $u_kB_2 u_k^* \cap B_1  \notin \mathcal{N}$.
\end{enumerate}
On the other hand, let  $\varepsilon>0$ be such that
$\pi(\mathcal{N}_{\varepsilon}(1_{B_1})) \subseteq \beta(\mathcal{N}) $.
Let $\{ e_k(s,i,j) \}_{1\leq s \leq l, 1\leq i ,j \leq k_s}$ denote  a system of  matrix units for $u_kB_2u_k^*\cap B_1$. Fix elements  $f_k(s,i,j)$ in $B_2$ such that
 $e_k(s,i,j)=u_kf_k(s,i,j)u_k^*$.
 Since $B_2$ is finite dimensional,
passing to a subsequence if necessary, we may assume that $\lim_k f_k(s,i,j)=f(s,i,j)$, for all $s,i$ and $j$. Using  property (1) of the sequence $(u_k)_{k\geq 1}$, we deduce 
 $$\lim_k e_k(s,i,j)  = \lim_k u_kf_k(s,i,j)u_{k}^*=u_0f(s,i,j)u_0^*  .$$
 Hence the element  $e(s,i,j)=u_0f(s,i,j)u^* $ belongs to $ u_0B_1u_0^* \cap B_1=B$.
 Use Lemma \ref{PerturbationDavidson}  and take $\delta_1 $ positive    such that whenever $C$ is a subalgebra  in $\SA(B_1)$  having a  system of matrix units $\{ e_C(s,i,j) \}_{s,i,j}$ satisfying
$\dist(e_C(s,i,j),B) < \delta_1$, for all $s,i$ and $j$,  then there is a unitary $Q$ in $\mathbb{U}(B_1)$ such that $\| Q-1_{B_1} \|< \varepsilon$ and $QCQ^* \subseteq B$. 
 Take   $k$  such that $\| e_{k}(s,i,j) - e(s,i,j)  \| < \delta_1$ for all $s,i$  and $j$. This implies $\dist(e_{k}(s,i,j),B) < \delta_1$ for all $s,i$ and $j$. We conclude there is a unitary $Q$ in $\mathbb{U}(B_1) $ such that $\| Q-1_{B_1} \|< \varepsilon$ and 
$Q^*(u_{k}B_2 u_{k}^*\cap B_1  )Q\subseteq B$. But  
$$\dim B =\dim u_{k}B_2u_{k}^* \cap B_1  =\dim Q^*(u_{k}B_2 u_{k}^*\cap B_1  )Q ,$$
where in the first equality we used that $u_k $ lies in $Z(B_1,B_2;[B]_{B_1})$.
Hence  $Q^*(u_{k}B_2 u_{k}^*\cap B_1  )Q=B$.  As a consequence,  
$$\beta(u_{k}B_2u_{k}^* \cap B_1)=\beta(QBQ^*)=\pi(Q) \in \beta(\mathcal{N}) .$$ 
 But the latter contradicts property (2) of $(u_k)_{k\geq 1}$.
\end{proof}

\begin{lemma}\label{CenterConti}
For $B$ in $\SA(B)$, the function $c:[B]_{B_1}\to [C(B)]_{B_1}$ given by $c(uBu^*)= uC(B)u^*$ is continuous. 
\end{lemma}
\begin{proof}
First, we must show the function $c$ is well defined. In other words we have to show $\Sb(B_1,B)\subseteq \Sb(B_1,C(B))$. But this follows directly  from the fact that any $u$ in $\Sb(B_1,B)$ defines a $*$--automorphism of $B$ and any $*$--automorphism leaves the center fixed.
Since $[B]_{B_1}$ and $[C(B)]_{B_1}$ are  homeomorphic to $\mathbb{U}(B_1)/\Sb(B_1,B)$ and $\mathbb{U}(B_1)/\Sb(B_1,C(B))$ respectively, it follows that $c$ is continuous if and only if  the function $\tilde{c}:\mathbb{U}(B_1)/\Sb(B_1,B) \to \mathbb{U}(B_1)/\Sb(B_1,C(B))$ given by $\tilde{c}(u\Sb(B_1,B))=u\Sb(B_1,C(B))$ is continuous. But the spaces $\mathbb{U}(B_1)/\Sb(B_1,B)$ and $\mathbb{U}(B_1)/\Sb(B_1,C(B))$  have the quotient topology induced by the canonical projections 
$$\pi_B: \mathbb{U}(B_1)\to \Sb(B_1,B) , \quad  \pi_{C(B)}:\mathbb{U}(B_1) \to \mathbb{U}(B_1)/\Sb(B_1,C(B)).$$
Thus $\tilde{c}$ is continuous if and only if $\pi_B \circ \tilde{c}$ is continuous. But $\pi_B\circ \tilde{c}=\pi_{C(B)}$, which is indeed continuous.
\end{proof}

We are  ready to find  local parameterizations  of $Z(B_1,B_2;[B]_{B_1})$.
\begin{prop}\label{EmbeddingALittlePiece}
Take $B $ a unital $C^*$-subalgebra in $B_1$ that is unitarily equivalent to a $C^*$-subalgebra of $B_2$. Fix an element $u_0$ in  $Z(B_1,B_2;[B]_{B_1})$. 
Then there is  a positive number $r$ and a continuous injective function 
$$\Psi: \mathcal{N}_r (u_0)\cap Z(B_1,B_2;[B]_{B_1}) \to \mathbb{R}^{d(C(B))}.$$
\end{prop}
\begin{proof}
Using that $Z(B_1,B_2;[B]_{B_1})=Z(B_1,B_2,[u_0B_2u_0^*\cap B_1]_{B_1})$, with no loss of generality we may  assume  $u_0B_2u_0^* \cap B_1 =B$.
Now, we use the manifold structure of  $[C(B)]_{B_1}$ and $Y(B_2;C(B))$ to construct $\Psi$. Note that if $Y(B_2,B)$ is nonempty then $Y(B_2,C(B))$ is nonempty as well. 
Let $d_1$ denote the dimension of $[C(B)]_{B_1}$ and let $d_2$ denote the dimension of the submanifold of $Y(B_2;C(B))$
that contains $u_0$.
Of course, we have $d_1+d_2\le d(C(B))$.
 
We use the local cross section result from previous section to parametrize $[C(B)]_{B_1}$. To ease notation  take $G=\mathbb{U}(B_1)$, $H=\Sb(B_1,C(B))$ and let $\pi$ denote the canonical quotient map from $G$ onto the left-cosets of $H$. By Proposition  \ref{LocalCrossSection} there are
 \begin{enumerate}
 \item $\mathcal{N}_G$, a compact neighborhood of $1$ in $G$,
 \item $\mathcal{N}_H$, a compact  neighborhood of $1$ in $H$,
 \item $\mathcal{N}_{G/H}$, a compact neighborhood of $\pi(1)$ in $G/H$,
 \item a continuous function $s:\mathcal{N}_{G/H} \to \mathcal{N}_G$ satisfying
 \begin{enumerate}
 \item $s(\pi(1))=1$ and $\pi(s(\pi(g)))=\pi(g)$ whenever $\pi(g)$ lies  in $\mathcal{N}_{G/H}$,
 \item the function
 \begin{eqnarray*}
 \mathcal{N}_{H} \times \mathcal{N}_{G/H}  & \to & \mathcal{N}_G, \\
 (h,\pi(g))  & \mapsto  & hs(\pi(g)),
 \end{eqnarray*}
 is an homeomorphism.
\end{enumerate}
\end{enumerate}
Since $G/H$ is a manifold of dimension $d_1$, we may assume there is a continuous injective map $\Psi_1: \mathcal{N}_{G/H} \to \mathbb{R}^{d_1}$.

Parametrizing  $Y(B_2;C(B))$ is easier. Since $u_0B_2u_0^*\cap B_1 =B$, $u_0$ belongs to  $Y(B_2;B)$. Take $r_1$ positive and a  diffeomorphism $\Psi_2$ from $Y(B_2;C(B))\cap \mathcal{N}_{r_1}(u_0)$ onto an open subset of $\mathbb{R}^{d_2}$.

Now that we have fixed parametrizations $\Psi_1$ and $\Psi_2$,  we can  parametrize  $Z(B_1,B_2;[B]_{B_1})$ around $u_0$.
Recall $[C(B)]_{B_1}$ has the topology induced by the bijection $\beta:[C(B)]_{B_1} \to G/H$, given by $\beta(uC(B)u^*)=\pi(u)$.
The function 
\begin{eqnarray*}
Z(B_1,B_2;[B]_{B_1} ) & \to  & [C(B)]_{B_1},\\
u  & \mapsto & c(uB_2u^*\cap B_1)
\end{eqnarray*}
is continuous by Lemma \ref{ContinuityIntersectionProp} and Lemma \ref{CenterConti}. Hence there is $\delta_2$ positive  such that $\beta(c(uB_2u^*\cap B_1))$ belongs to $\mathcal{N}_{G/H}$, whenever $u$ lies in the intersection $Z(B_1,B_2;[B]_{B_1}) \cap \mathcal{N}_{\delta_2}(u_0)$.
For a unitary $u$ in $Z(B_1,B_2;[B]_{B_1}) \cap \mathcal{N}_{\delta_2}(u_0)$ define
$$
q(u):=s(\beta(c(uB_2u^*\cap B_1))).
$$
We note that $q(u_0)=1$, $q(u)$ lies in $G$ and that the map $u \mapsto q(u)$ is continuous. The main property of $q(u)$ is that 
\begin{equation}\label{MainPropertyq(u)}
c(uB_2u^*\cap B_1) =q(u)c(B)q(u)^*.
 \end{equation}
  Indeed, for $u$ in $Z(B_1,B_2;[B]_1)\cap \mathcal{N}_{\delta_2}(u_0)$ there is  a unitary  $v$ in $G$ with the property  $uB_2u^*\cap B_1 = vBv^*$. Hence $c(uB_2 \cap B_1)=vC(B)v^*$. Since $\|  u  - u_0 \| <\delta_2$, $\beta(c(uB_2u^*\cap B_1)) $ lies in $\mathcal{N}_{G/H}$. Hence  $\beta(c(uB_2u^* \cap B_1))=\pi(v)$ lies in $\mathcal{N}_{G/H}$. Using  the fact that $s$ is a local section on $\mathcal{N}_{G/H}$ (property  (4a) above) we deduce  $\pi(s(\pi(v)))=\pi(v)$ .
  
 On the other hand, by definition of $q(u)$ we have
 $$
 \pi(s(\pi(v)))=\pi(s(\beta(uB_2u^*\cap B_1)))=\pi(q(u)).
 $$
As a consequence, $\pi(v)=\pi(q(u))$ i.e. $v^*q(u) $ belongs to $\Sb(B_1,B)$ which is just another way to say (\ref{MainPropertyq(u)}) holds.
At last we are ready to  find $r$. Continuity of the map $u \mapsto q(u)$ gives a positive $\delta_3$, less that $\delta_2$, such that $\| q(u)-1 \|< \frac{\delta_1}{2}$  whenever $u$ lies  in $Z(B_1,B_2;[B]_{B_1}) \cap \mathcal{N}_{\delta_3}(u_0)$. Define $r=\min\{  \frac{\delta_1}{2}, \delta_3 \}$.
The first thing we notice is that $q(u)^*u$ belongs to  $Y(B_2;C(B))\cap \mathcal{N}_{\delta_1}(u_0)$ whenever $u$  lies in $Z(B_1,B_2;[B]_{B_1})\cap \mathcal{N}_{\delta}(u_0)$.  Indeed, from
$$
q(u)c(B)q(u)^*=c(uB_2u^*\cap B_1) \subseteq uB_2 u^*
$$
we obtain  $q(u)^*u \in Y(B_2;c(B))$ and  a standard computation, using $\| q(u)-1 \|< \frac{\delta_1}{2}$, shows $\| q(u)^*u-u_0 \| < \delta_1$.
Hence we are allowed to take $\Psi_2(q(u)^*u)$. 
Lastly,  for $u $ in $Z(B_1,B_2;[B]_{B_1}) \cap \mathcal{N}_{\delta}(u_0)$ define 
$$
\Psi(u):=( \Psi_1( \beta(c(uB_2u^* \cap B_1 )) ), \Psi_2(q(u)u^* )).
$$
It is clear  that $\Psi$ is continuous.

Now we show $\Psi$ is injective. If $\Psi(u_1)=\Psi(u_2)$, for  two element  $u_1$ and  $u_2$ in  $Z(B_1,B_2;[B]_{B_1})$, then
\begin{eqnarray}
\Psi_1(\beta(c(u_1B_2u_1^* \cap B_1))) &=& \Psi_1(\beta(c(u_2B_2u_2^*\cap B_1))), \label{EqualityAlgebras} \\
\Psi_2(q(u_1)u_1^*) &=& \Psi_2(q(u_2)u_2^*) . \label{EqualityTwistingUnitary}
\end{eqnarray}
From (\ref{EqualityAlgebras}) and definition of $q(u)$ it follows that $q(u_1)=q(u_2)$ and from equation (\ref{EqualityTwistingUnitary}) we conclude $u_1=u_2$.
\end{proof}

\begin{prop}\label{EmbeddingALittlePieceCase2}
Take $B $ a unital $C^*$-subalgebra of  $B_1$ such that it is unitarily equivalent to a $C^*$-subalgebra of $B_2$.  Fix an element $u_0$ in  $Z(B_1,B_2;[B]_{B_1})$. 

There is  a positive number $r$ and a continuous injective function 
$$\Psi: \mathcal{N}_r (u_0)\cap Z(B_1,B_2;[B]_{B_1}) \to \mathbb{R}^{d(B)}$$.
\end{prop}

The proof of Proposition \ref{EmbeddingALittlePieceCase2} is similar to that of  Proposition \ref{EmbeddingALittlePiece}, so we omit it.

We now begin showing density in $\mathbb{U}(M_N)$ of certain sets of unitaries.
\begin{lemma}\label{DensityZComplementSimpleCase}
Assume $B_1$ and $B_2$ are simple. If $B\not= \mathbb{C}$ is a unital $C^*$-subalgebra  of $B_1$  and it is unitarily equivalent to a $C^*$-subalgebra of $B_2$  then $Z(B_1,B_2;[B]_{B_1})^c$ is dense.
\end{lemma}
\begin{proof}
Firstly we notice that $\dim \mathbb{U}(B_1)+\dim \mathbb{U}(B_2)<N^2$. Indeed, if $B_i$ is $*$--isomorphic to $M_{k_i}$, $i=1,2$ and $m_i=\mu(M_N,B_i)$ then $\dim \mathbb{U}(B_1)+ \dim \mathbb{U}(B_2) =N^2(1/m_2^2+1/m_2^2)<N^2$.
Secondly we will prove that for any $u$ in $Z(B_1,B_2;[B]_{B_1})$ there is a natural number $d_u$, with $d_u<N^2$, a positive number $r_u$ and a continuous injective function $\Psi_u:\mathcal{N}_{r_u}(u)\cap Z(B_1,B_2;[B]_{B_1}) \to \mathbb{R}^{d_u}$.
We will consider two cases.

\medskip
\noindent
{\em Case (1):} $B$ is not simple. 
Take $d_u=d(C(B))$. Since $C(B)\not= \mathbb{C}$,  Proposition \ref{dBLessN2SimpleCase} implies  $d(C(B))< N^2$. Take $r_u$ and $\Psi_u$ as required to exist by Proposition \ref{EmbeddingALittlePiece}.

\medskip
\noindent
{\em Case (2):} $B$ is simple.
Take $d_u=d(B)$. Since $B\not= \mathbb{C}$, $B$  contains a unital  $C^*$-subalgebra isomorphic to $\mathbb{C}^2$,  call it $C$. Lemma \ref{BSimpleCIsomC2} implies $d(B) \leq d(C)$ and  Lemma \ref{dBLessN2SimpleCase} implies  $d(C)< N^2$. Take $r_u$ and $\Psi_u$ the positive number and continuous injective function from Proposition \ref{EmbeddingALittlePieceCase2}.

We will show that $U \cap Z(B_1,B_2;[B]_{B_1})^c\not= \emptyset $, for any nonempty open subset $U \subseteq \mathbb{U}(M_N)$. First notice that if  the intersection  $U \cap  (  \bigcup_{u\in Z(B_1,B_2;[B]_{B_1}) } \mathcal{N}_{r_u}(u) )^c$ is nonempty then 
we are done. Thus we may assume $ U \subseteq  \bigcup_{u\in Z(B_1,B_2;[B]_{B_1}) } \mathcal{N}_{r_u}(u) $. Furthermore, by making $U$ smaller, if necessary,  we may assume there is $u$ in $Z(B_1,B_2;[B]_{B_1})$ such that $U \subseteq \mathcal{N}_{r_u}(u)$.

For sake of contradiction assume $U \subseteq Z(B_1,B_2;[B]_{B_1})$. We may take an open subset $V$, contained in  $U$, small enough so that $V$ is diffeomorphic to an open connected set $\mathcal{O}$ of $\mathbb{R}^{N^2}$. Let $\varphi : \mathcal{O} \to V$ be a diffeomorphism.
It follows we have a continuous injective function
$$
\xymatrix{
 \mathbb{R}^{N^2} \supseteq \mathcal{O}   \ar[r]^-\varphi   & V  \ar[r]^{\Psi_u} &  \mathbb{R}^{d_u} \ar@{^{(}->}[r] & \mathbb{R}^{N^2} 
}
.$$

By the Invariance of Domain Theorem, the image of this map must be open in $\mathbb{R}^{N^2}$. But  this is a  contradiction since the image is contained in $\mathbb{R}^{d_u}$ and $d_u <N^2$. We conclude $U\cap Z(B_1,B_2;[B]_{B_1})^c \not= \emptyset$.
\end{proof}

\begin{lemma}\label{DensityZComplementBigCenterCase}
Suppose $\dim C(B_1) \geq 2$  and $B_1$ is $*$--isomorphic to
$$
M_{N/\dim C(B_1)}\oplus \cdots \oplus M_{N/\dim C(B_1)}.
$$
Assume one of the following cases holds:
\begin{enumerate}
\item $\dim C(B_2)=1$,
\item $B_1$ is $*$--isomorphic to
$$
M_{N/2} \oplus M_{N/2}
$$
and $B_2$ is $*$--isomorphic to 
$$
M_{N/2}\oplus M_{N/(2k)},
$$
where $k \geq 2$.
\item $\dim C(B_2)\geq 3$  and $B_2$ is $*$--isomorphic to
$$
M_{N/\dim C(B_2)}\oplus \cdots \oplus M_{N/\dim C(B_2)}.
$$
\end{enumerate}
 Then for any  $B\not= \mathbb{C}$ unital $C^*$-subalgebra of  $B_1$  such that it is unitarily equivalent to a $C^*$-subalgebra  of $B_2$, $Z(B_1,B_2;[B]_{B_1})^c$ is dense.
\end{lemma}
\begin{proof}
The proof of  Lemma \ref{DensityZComplementBigCenterCase} is exactly as the proof  of \ref{DensityZComplementSimpleCase} 
but using Lemma \ref{dBLessN2NonSimpleCase} instead of Lemma \ref{dBLessN2SimpleCase} .
\end{proof}

At this point if the sets $Z(B_1,B_2;[B]_{B_1})$ were closed one could  conclude immediately that $\Delta(B_1,B_2)$ is dense. Unfortunately they may not be closed. What saves the day is the fact that  we can control the closure of $Z(B_1,B_2;[B]_{B_1})$ with sets of the same form i.e. sets like $Z(B_1,B_2;[C]_{B_1})$ for a suitable  finite family of  subalgebras $C$. We make this statement clearer with the definition of an order on $\SA(B_1)$.
\begin{definition}
On $\SA(B_1)/\sim_{B_1}$ we define a partial order as follows:
$$
[B]_{B_1} \leq [C]_{B_1} \Leftrightarrow \exists D \in \SA(C) : D\sim_{B_1} B .
$$
\end{definition}

\begin{prop}\label{ClosureZContainedUnionOfZ}
For any $B$  in $\SA(B_1)$,
$$
\overline{Z(B_1,B_2;[B]_{B_1})} \subseteq \bigcup_{[C]_{B_1} \geq [B]_{B_1}} Z(B_1,B_2;[C]_{B_1}).
$$
\end{prop}
\begin{proof}
Let $(u_k)_{k\geq 1}$ be a sequence in $Z(B_1,B_2;[B]_{B_1})$ and $u$ in $\mathbb{U}(M_N)$ such that $\lim_k\|u_k-u\|=0$.
Pick $q_k$ in $\mathbb{U}(M_N)$ such that $q_kBq_k^*=u_kB_2u_k^*\cap B_1$.
Let $\{ f_k(s,i,j)\}_{s,i,j}$ be a matrix unit for $u_kB_2u_k^*\cap B_1$ and take elements $e_k(s,i,j)$ in $B_2$ such that
$f_k(s,i,j)=u_ke_k(s,i,j)u_k^*$. Since $B_2$ is finite dimensional, passing to a subsequence if necessary,
we may assume
$\lim_k  f_{k}(s,i,j)=f(s,i,j)\in B_2$ and $\lim_k u_{k}e_{k}(s,i,j)u_{k}^*=ue(s,i,j)u^*$ for some $e(s,i,j)\in B_1$,
for all $s$, $i$ and $j$. 
It follows that $\lim_k \dist(f_{k}(s,i,j), uB_2u^*\cap B_1)=0$.
Hence, from Lemma \ref{PerturbationDavidson},  for  large $k$, there is $q$ in $\mathbb{U}(M_N)$ so that
$q(u_{k}B_2u_{k}^*\cap B_1)q^*=qq_{k}Bq_{k}^*q^*$ is contained in $uB_2u^*\cap B_1$.
We conclude $[uB_2u^*\cap B_1]_{B_1}\geq [B]_{B_1}$ and since $u$ lies in $Z(B_1,B_2;[uB_2u^*\cap B_1])$ the proof is complete.
\end{proof}

\begin{lemma}\label{NonDensityZAndChains}
Assume one of the following conditions holds:
\begin{enumerate}
\item $\dim C(B_1)=1=\dim C(B_2)$,
\item $\dim C(B_1)  \geq 2$, $\dim C(B_2)=1$ and $B_1$ is $*$--isomorphic to
$$
M_{N/\dim C(B_1)}\oplus \cdots \oplus M_{N/\dim C(B_1)},
$$
\item $\dim C(B_1) = 2=\dim C(B_2)$, $B_1$ is $*$--isomorphic to
$$
M_{N/2} \oplus M_{N/2},
$$
and $B_2$ is $*$--isomorphic to
$$
M_{N/2}\oplus M_{N/(2k)}
$$
where $k \geq 2$,
\item $\dim C(B_1)\geq 2$, $\dim C(B_2) \geq 3$  and, for $i=1,2$, $B_i$ is $*$--isomorphic to
$$
M_{N/\dim C(B_i)}\oplus \cdots \oplus M_{N/ \dim C(B_i)}.
$$
\end{enumerate} 
Take $B$ a unital $C^*$-subalgebra  of $B_1$ such that it is unitarily equivalent to a  $C^*$-subalgebra of $B_2$ . 
If  $\quad \overline{Z(B_1,B_2;[B]_{B_1})}^c$ is not dense and $B\not= \mathbb{C}$ then there is a subalgebra  $C$  in $\SA(B_1)$ such that
$[C]_{B_1}> [B]_{B_1}$ and $\quad \overline{Z(B_1,B_2;[C]_{B_1})}^c$ is not dense.
\end{lemma}
\begin{proof}
We proceed by contrapositive. Thus,  assume $\overline{Z(B_1,B_2;[C]_{B_1})}^c$ is dense  for all $[C]_{B_1}>[B]_{B_1}$.
Since the set $\{ [C]_{B_1}: [C]_{B_1}>[B]_{B_1} \}$ is finite,
$$\bigcap_{[C]_{B_1}>[B]_{B_1}}\overline{Z(B_1,B_2;[C]_{B_1})}^c$$
 is open and dense.  
Furthermore, Lemma \ref{DensityZComplementSimpleCase} or Lemma \ref{DensityZComplementBigCenterCase} implies $Z(B_1,B_2;[B]_{B_1})^c$ is dense. Hence the intersection
 $$
 Z(B_1,B_2;[B]_{B_1})^c \cap \bigcap_{[C]_{B_1}>[B]_{B_1}} \overline{Z(B_1,B_2;[C]_{B_1})}^c
 $$
is dense. But this along with Proposition \ref{ClosureZContainedUnionOfZ} implies $\overline{Z(B_1,B_2;[B]_{B_1})}^c$ is dense.
\end{proof}

\begin{lemma}\label{DensityZClosureComplement}
Assume one of the conditions (1)--(4) of Lemma~\ref{NonDensityZAndChains} holds.
Then for any $B\not= \mathbb{C}$,  unital $C^*$-subalgebra of $B_1$ that is unitarily equivalent to a $C^*$-subalgebra of  $B_2$, the set  $\overline{Z(B_1,B_2;[B]_{B_1})}^c$ is dense.
\end{lemma}
\begin{proof}
Assume $\overline{Z(B_1,B_2;[B]_{B_1})}^c$ is not  dense. By Lemma \ref{NonDensityZAndChains} there is $[C]_{B_1}>[B]_{B_1}$ such that
$\overline{Z(B_1,B_2;[C]_{B_1})}^c$ is not dense. We notice that again we are in the same condition to apply   Lemma \ref{NonDensityZAndChains}, since $[C]_{B_1}>[B]_{B_1}>[\mathbb{C}]_{B_1}$. In this way we can construct chains, in $\SA(B_1)/\sim_{B_1}$, of length arbitrarily large, but this can not be since it is finite. 
\end{proof}

At last we  can give a proof of  Theorem \ref{DensityOfSmallIntersection}.

\begin{proof}[Proof of Theorem \ref{DensityOfSmallIntersection}]
A direct computation shows that
$$
\Delta(B_1,B_2)=\bigcap_{[B]_{B_1}  > [\mathbb{C}]_{B_1}} Z(B_1, B_2, [B]_{B_1})^c .
$$

Thus
$$
\Delta(B_1,B_2) \supseteq \bigcap_{[B]_{B_1}  > [\mathbb{C}]_{B_1}} \overline{Z(B_1, B_2, [B]_{B_1})}^c .
$$
Now, by  Lemma  \ref{DensityZClosureComplement},  whenever   $[B]_{B_1} > [\mathbb{C}]_{B_1} $, the set $\overline{Z(B_1,B_2,[B]_{B_1})}^c$ is dense. Hence  $\Delta(B_1,B_2)$ is dense.
\end{proof}

\section{Primitivity }\label{PrimitivitySection}

During this section, unless stated otherwise,  $A_1\not= \mathbb{C}$ and $A_2\not= \mathbb{C}$ denote  two nontrivial,  separable, residually finite dimensional $C^*$-algebras. Our goal  is to prove $A_1*A_2$ is primitive, except for the case $A_1=\mathbb{C}^2=A_2$. Two main ingredients are used. Firstly, the perturbation results from previous section. Secondly, the fact that $A_1*A_2$ has a   separating family of finite dimensional $*$--representations,  a result due to Excel and Loring, ~\cite{Exel&Loring}. 

Before we start proving results about primitivity, we want to consider the case $\mathbb{C}^2*\mathbb{C}^2$. This is a well studied $C^*$-algebra;  see for instance  ~\cite{Blackadar},   ~\cite{Pedersen} and ~\cite{Raeburn&Sinclair}. It is known  that $\mathbb{C}^2*\mathbb{C}^2$ is $*$--isomorphic to the $C^*$-algebra of continuous $M_2$-valued functions on the closed interval $[0,1]$,  whose values at 0 and 1 are diagonal matrices. As a consequences its center is not trivial. Since the center of any primitive $C^*$-algebra is trivial, we conclude $\mathbb{C}^2*\mathbb{C}^2$ is not primitive.

\begin{definition}
We denote by $\iota_j$  the inclusion $*$--homomorphism from $A_j$ into $A_1*A_2$.
Given a unital $*$--representation $\pi:A_1*A_2 \to \mathbb{B}(H)$, we define $\pi^{(1)}=\pi \circ \iota_1$ and $\pi^{(2)}=\pi \circ \iota_2$. Thus, with this notation, we have $\pi=\pi^{(1)}* \pi^{(2)}$.  For a unitary $u$ in $\mathbb{U}(H)$ we call the $*$--representation $\pi^{(1)}*(\Ad u\circ \pi^{(2)})$, a perturbation of $\pi$ by $u$.
\end{definition}

\begin{remark}\label{Perturbations&Irreducibility}
The $*$--representation $\pi^{(1)}*(\Ad u \circ \pi^{(2)} )$ is irreducible if and  only if
$$
u\pi^{(2)}(A_2)'u^* \cap \pi^{(1)}(A_1)' =\mathbb{C}.
$$
where $(\pi^{(1)}(A_1))'$ denotes de commutant of $\pi^{(1)}(A_1) $ in $\mathbb{B}(H)$.
\end{remark}

Our first goal is to perturb a given finite dimensional $*$--representation of $A_1*A_2$ into an irreducible one. Of course, 
the example $\mathbb{C}^2*\mathbb{C}^2$ shows that in general
 this can't be done so we have to find conditions that guarantee it. We start with the case $A_1$ and $A_2$ finite
 dimensional and later, built on the finite dimensional case, we continue with the residually finite dimensional case. For the finite dimensional case crucial information is given by the ranks of minimal central projections on $A_1$ and $A_2$.

\begin{definition}\label{def:RCP}
Assume $A_1$ and $A_2$ are finite dimensional and
let $\rho:A_1*A_2\to\mathbb{B}(H)$ be a unital, finite dimensional representation.
We say that $\rho$ satisfies the {\em Rank of Central Projections condition} (or {\em RCP condition}) if
for both $i=1,2$, the rank of $\rho(p)$ is the same for all minimal projections $p$ of the center $C(A_i)$ of $A_i$,
(but they need not agree for different values of $i$).
\end{definition}

The RCP condition for $\rho$, of course, is really about the pair of representations $(\rho^{(1)},\rho^{(2)})$.
However, it will be convenient to express it in terms of $A_1*A_2$.
In any case, the following two lemmas are clear.

\begin{lemma}\label{lem:RCPunitary}
Suppose $A_1$ and $A_2$ are finite dimensional,
$\rho:A_1*A_2\to\mathbb{B}(H)$ is a finite dimensional representation
that satisfies the RCP condition and $u\in\mathbb{U}(H)$.
Then the representation $\rho^{(1)}*(\Ad u\circ\rho^{(2)})$ of $A_1*A_2$
also satisfies the RCP condition.
\end{lemma}

\begin{lemma}\label{lem:RCPsum}
Suppose $A_1$ and $A_2$ are finite dimensional,
$\rho:A_1*A_2\to\mathbb{B}(H)$ and $\sigma:A_1*A_2\to\mathbb{B}(K)$
are finite dimensional representations that satisfy the RCP condition.
Then $\rho\oplus\sigma:A_1*A_2\to\mathbb{B}(H\oplus K)$ also satisfies the RCP condition.
\end{lemma}

The following is clear from Lemma \ref{PartialMultiplicities}.
\begin{lemma}\label{LemmaMultiplicityCenter}
Assume $A$ is a finite dimensional $C^*$-algebra $*$--isomorphic to
$\bigoplus_{j=1}^l M_{n(j)}$ 
and  take $\pi:A\to \mathbb{B}(H)$  a unital  finite dimensional $*$--representation. 
Let $\mu(\pi)=[m(1),\ldots,m(l)]$ and let $\tilde{\pi}$ be the restriction of $\pi$ to the center of $A$. Then
$$
\mu(\tilde{\pi})=[m(1)n(1),\dots ,m(l)n(l) ].
$$
\end{lemma}

The next lemma will help us to prove that the RCP condition is easy to get. 

\begin{lemma}\label{AddingPiecesToFinDimRep}
Assume $A$ is a finite dimensional $C^*$-algebra and $\pi:A \to \mathbb{B}(H)$ is a  unital
finite dimensional $*$--representation. Let
$$
\mu(\pi)=[m(1),\dots  ,m(l)].
$$
For any nonnegative integers $q(1),\dots, q(l)$ there is a finite dimensional unital $*$--representation $\rho:A \to \mathbb{B}(K)$ such that
$$
\mu(\pi \oplus \rho)=[m(1)+q(1),\dots ,m(l)+q(l)].
$$
\end{lemma}
\begin{proof}
Write $A$ as
$$
A=\bigoplus_{i=1}^l A(i) 
$$
where $A(i)=\mathbb{B}(V_i)$ for $V_i$ finite dimensional.
For $1\leq i \leq l$, let $p_{i}:A\to A(i)$ denote the canonical projection onto $A(i)$. Notice that $p_i$ is a unital $*$--representation of $A$.
Define
$$
\rho:= \bigoplus_{i=1}^l \underbrace{( p_i \oplus \cdots \oplus p_i )}_{q(i)-\textrm{times}}:A\to \bigoplus_{i=1}^l A(i)^{q(i)}\subseteq\mathbb{B}(K),
$$
where $K=\bigoplus_{i=1}^l(V_i^{\oplus q_i})$.
Then $\rho$ is a unital $*$--representation of $A$ on $K$ and 
$$
\mu(\pi\oplus \rho)=[m(1)+q(1),\dots , m(l)+q(l)].
$$
\end{proof}

The next lemma takes slightly more work and is essential to our construction.

\begin{lemma}\label{lem:RCP}
Assume $A_1$ and $A_2$ are finite dimensional. Given a unital finite dimensional $*$--representation $\pi:A_1*A_2 \to \mathbb{B}(H)$, there is a finite dimensional Hilbert space $\hat{H}$ and   a unital $*$--representation 
$$\hat{\pi}:A_1*A_2 \to \mathbb{B}(\hat{H})$$
such that $\pi\oplus \hat{\pi}$ satisfies the RCP condition.
\end{lemma}
\begin{proof}
For $i=1,2$, let $l_i=\dim C(A_i)$, let $A_i$ be $*$--isomorphic to 
$\bigoplus_{j=1}^{l_i}M_{n_i(j)}$
and  write
\[
\mu(\pi^{(i)})=[ m_i(1),\dots , m_i(l_i)  ].
\]
Take
$n_i=\lcm(n_i(1),\dots, n_i(l_i))$
and integers $r_i(j)$, such that $r_i(j)n_i(j)=n_i$,  for $1\leq j  \leq l_i$. Take a positive integer $s$ such that $sr_i(j) \geq m_i(j)$ for all $i=1,2$ and $1\leq j \leq l_i$.
Use Lemma \ref{AddingPiecesToFinDimRep} to find a unital  finite dimensional $*$--representation $\rho_i:A_i \to \mathbb{B}(K_i)$, $i=1,2$ such that
$$
\mu(\pi^{(i)}\oplus \rho_i)=[sr_i(1),\dots , sr_i(l_i)].
$$
Letting $\kappa_i$ denote the restriction of $\pi^{(i)}\oplus \rho_i$ to $C(A_i)$, from Lemma \ref{LemmaMultiplicityCenter}
we have
$$
\mu(\kappa_i)=[sr_i(1)n_i(1),\dots , sr_i(l_i)n_i(l_i)]=[sn_i,sn_i,\ldots,sn_i].
$$

The $*$--representations $(\pi^{(1)}\oplus \rho_1)$ and $(\pi^{(2)} \oplus \rho_2)$ are almost what we want,
but they may take values in Hilbert spaces with different dimensions.
To take care of this, we take multiples of them.
Let $N=\lcm(\dim(H\oplus K_1), \dim(H\oplus K_2))$, find positive  integers $k_1$ and $k_2$ such that
$$
N=k_1 \dim(H \oplus K_1)=k_2\dim(H\oplus K_2)
$$
and consider the Hilbert spaces $(H\oplus K_i)^{\oplus k_i}$, whose dimensions agree for
$i=1,2$. 
Then
$$
\dim( K_1\oplus (H\oplus K_1)^{\oplus(k_1-1)} )=\dim( K_2\oplus (H\oplus K_2)^{\oplus(k_2-1)} )
$$
and there is a unitary operator  
$$
U: K_2\oplus (H\oplus K_2)^{\oplus(k_2-1)} \to K_1\oplus (H\oplus K_1)^{\oplus(k_1-1)}.
$$ 
Take 
\begin{eqnarray*}
\hat{H}&:=&K_1\oplus (H+K_1)^{\oplus(k_1-1)},\\
\hat{\pi}_1&:=&\rho_1 \oplus ( \pi^{(1)} \oplus \rho)^{\oplus(k_1-1)},\\
\sigma_1&:=&\pi^{(1)}\oplus \hat{\pi}_1,\\
\hat{\pi}_2&:=& \Ad U \circ (\rho_2 \oplus ( \pi^{(2)} \oplus \rho )^{\oplus(k_2-1)}),\\
\sigma_2&:=&\pi^{(2)}\oplus \hat{\pi}_2,\\
 \hat{\pi}&:=&\hat{\pi}_1* \hat{\pi}_2 .
\end{eqnarray*}
Then $\sigma_1*\sigma_2=(\pi^{(1)}\oplus \hat{\pi}_1)*(\pi^{(2)}  \oplus \hat{\pi}_2)=\pi\oplus \hat{\pi}$.
We have $\mu(\sigma_i)=[k_isr_i(1), \dots , k_isr_i(l_i)]$.
Let $\tilde{\sigma}_i$ denote the restriction of $\sigma_i$  to $C(A_i)$.
From Lemma \ref{LemmaMultiplicityCenter} we have
$$
\mu(\tilde{\sigma}_i)=[k_isr_i(1)n_i(1), \dots , k_isr_i(l_i)n_i(l_i)]=[k_isn_i,\dots , k_isn_i].
$$
\end{proof}

The purpose of the next definition and lemma is to emphasize an important property about 
$*$--representations satisfying the RCP.

\begin{definition}\label{def:DPI}
A $*$--representation $\pi:A_1*A_2 \to \mathbb{B}(H)$ is said to be
{\em densely perturbable to an irreducible $*$-representation}, abbreviated
{\em DPI},
if the set
\[
\Delta(\pi):=\{u \in \mathbb{U}(H): \pi^{(1)}(A_1)'\cap ( u \pi^{(2)}(A_2)'u^* ) =\mathbb{C}\}
\]
is norm dense in $\mathbb{U}(H)$. Here the commutants are taken with respect to
$\mathbb{B}(H)$.
\end{definition}

The next lemma shows that any $*$--representation satisfying the RCP is DPI.

\begin{lemma}\label{lem:RCPimpliesDPI}
Assume $A_1$ and $A_2$ are finite dimensional C*-algebras and $(\dim(A_1)-1)(\dim(A_2)-1)\geq 2$.
If $\rho:A_1*A_2 \to \mathbb{B}(H)$, with $H$ finite dimensional, satisfies the Rank of Central Projections condition,
then $\rho$ is DPI.
\end{lemma}

\begin{proof}
Since $(\dim(A_1)-1)(\dim(A_2)-1)\geq 2$, and after interchanging $A_1$ and $A_2$, if necessary, 
one of the following must hold:
\begin{enumerate}[(1)]
\item $A_1$ and $A_2$ are simple,
\item $\dim C(A_1) \geq 2$ and $A_2$ is simple,
\item for $i=1,2$, $A_i=M_{n_i(1)}\oplus M_{n_i(2)}$, with $n_2(2) \geq 2$,
\item $\dim C(A_1) \geq 2$,  $\dim C(A_2) \geq 3$.
\end{enumerate}

In case~(1), take $B_i=\rho^{(i)}(A_i)'$, $i=1,2$. 

In case~(2), let $B_1=\rho^{(1)}(C(A_1))'$ and $B_2=\rho^{(2)}(A_2)'$.
Notice that $\dim C(B_2)=1$, $\dim C(B_1)=\dim C(A_1) \geq 2$
and, by the RCP assumption, $B_1$ is $*$--isomorphic to 
$M_{\dim H /\dim C(B_1)}\oplus \cdots \oplus M_{\dim H  /\dim C(B_1)}.$

In case~(3),
let $B_1=\rho^{(1)}(C(A_1))'$ and $B_2=\rho^{(2)}(\mathbb{C}\oplus M_{n_2(2)})'$.
By the RCP assumption, $B_1$ is $*$--isomorphic to
$$
M_{\dim H /2} \oplus M_{\dim H /2} 
$$
and   $B_2$ is $*$--isomorphic to
$$
M_{\dim H /2}\oplus M_{\dim H /(2n_2(2))}.
$$

In case~(4),
let $B_i=\rho^{(i)}(C(A_i))'$ for $i=1,2$.
Then $\dim C(B_1)=\dim C(A_1) \geq 2$, $\dim C(B_2)=\dim C(A_2) \geq 3$ and, for $i=1,2$, RCP implies
$B_i$ is $*$--isomorphic to
$$
M_{\dim H /\dim C(B_i)}\oplus \cdots \oplus M_{\dim H /\dim C(B_i)}.
$$

Now define 
\[
\Delta(B_1,B_2):=\{u\in \mathbb{U}(H): B_1\cap \Ad u (B_2) =\mathbb{C}  \}.
\]
and notice that in all four cases $\Delta(B_1,B_2) \subseteq\Delta( \rho)$.
By Theorem \ref{DensityOfSmallIntersection}, the 
set $\Delta(B_1,B_2)$ is dense in all the four cases.
\end{proof}

A downside of the  DPI property is that it is not stable under direct sums.
However, it is stable under perturbations.

\begin{remark}\label{Remark:DPI}
If $\pi:A_1*A_2 \to \mathbb{B}(H)$ is DPI, then for any $u$ in $\mathbb{U}(H)$,
$\pi^{(1)}* ( \Ad u \circ \pi^{(2)})$ is also DPI. Indeed, this follows from the 
identity
\[
\Delta\big(\pi^{(1)}*  (\Ad u \circ \pi^{(2)})\big)=\Delta(\pi)u^*.
\]
\end{remark}

From Lemma \ref{lem:RCP} we obtain the following.

\begin{lemma}\label{Lemma:Perturbation}
For any unital finite dimensional $*$-representation $\pi:A_1*A_2 \to \mathbb{B}(H)$, 
there is a unital finite dimensional $*$-representation $\hat{\pi}:A_1*A_2 \to \mathbb{B}(\hat{H})$
such that $\pi \oplus \hat{\pi}$ is DPI.
\end{lemma}
\begin{proof}
The assumption $(\dim(A_1)-1)(\dim(A_2)-1)\geq 2$ implies there is a
unital finite dimensional $*$-representation $\vartheta:A_1*A_2\to \mathbb{B}(H_0)$,
such that $(\dim(\vartheta^{(1)}(A_1)) -1)(\dim(\vartheta^{(2)}(A_2))-1) \geq 2 $.
Consider the unital C*-subalgebras of $\mathbb{B}(H\oplus H_0)$,
$D_i=(\pi\oplus \vartheta)^{(i)}(A_i)$, $i=1,2$, and notice that $(\dim(D_1) -1)(\dim(D_2)-1) \geq 2$. 
Let  $\theta:D_1*D_2\to \mathbb{B}(H \oplus H_0 ) $ be the unital $*$-representation
induced by the universal property of $D_1*D_2$ via the unital inclusions
$D_i \subseteq \mathbb{B}(H\oplus H_0)$. 
Lemma \ref{lem:RCP} implies  there is a unital finite 
dimensional $*$-representation
$\rho:D_1*D_2 \to \mathbb{B}(K)$ such that $\theta \oplus \rho$  satisfies the RCP condition, 
so by Lemma \ref{lem:RCPimpliesDPI} is DPI.

Let $j_i: D_i \to D_1*D_2$, $i=1,2$, be the inclusion $*$-homomorphism from the
definition of unital full free product. Now consider the unital $*$-homomorphism
$\sigma=(j_1\circ (\pi \oplus \vartheta)^{(1)})*(j_2 \circ (\pi\oplus \vartheta
)^{(2)}):A_1*A_2 \to D_1*D_2$.
Now just take  $\hat{H}=H_0\oplus K$ and $\hat{\pi}=\vartheta \oplus (\rho  \circ \sigma)$. 
In order to show $\pi \oplus \hat{\pi}$ is DPI we just need to show that, for $i=1,2$,
$(\pi \oplus \hat{\pi})^{(i)}(A_i)=(\theta \oplus \rho )^{(i)}(D_i),$
but this is a direct computation.
\end{proof}

The proof of next lemma is a standard approximation argument and we omit it.

\begin{prop}\label{Perturbations&ApproximationOnGivenFiniteSets}
Let $A_1$ and $A_2$ be two unital $C^*$-algebras. Given  a non zero element $x$ in $A_1*A_2$ 
and a positive number $\varepsilon$, there is a positive number  $\delta=\delta(x,\varepsilon)$ such 
that for any   $u$ and $v$ in $\mathbb{U}(H)$ satisfying $\| u-v \|< \delta$ and any 
unital $*$-representations $\pi:A_1*A_2 \to \mathbb{B}(H)$, we have
$$
\|  (\pi^{(1)}*    ( \Ad v\circ \pi^{(2)})   )(x)- (\pi^{(1)}*    ( \Ad u\circ \pi^{(2)})   )(x)  \|< \varepsilon.
$$
\end{prop}

Here is our main theorem.

\begin{thm}\label{Primitivity}
Assume $A_1$ and $A_2$ are unital, separable, residually finite dimensional $C^*$-algebras with
$(\dim(A_1)-1)(\dim(A_2)-1) \geq 2$.
Then $A_1*A_2$ is primitive.
\end{thm}
\begin{proof}
By the result of Exel and Loring in~\cite{Exel&Loring}, there is a separating sequence  
$(\pi_j:A_1*A_2 \to \mathbb{B}(H_j))_{j\geq 1} $, of finite dimensional unital $*$-representations. 
For later use in constructing an essential representation of $A_1*A_2$, i.e., a $*$-representation
with the property that zero is the only
compact operator in its image, we modify $(\pi_j)_{j \geq 1}$, if necessary, so that  that each 
$*$-representation is repeated infinitely many times.
 
By recursion and using Lemma \ref{Lemma:Perturbation}, we define a 
sequence 
\[
\hat{\pi}_j:A_1*A_2 \to \mathbb{B}(\hat{H}_j),\quad(j\ge1)
\]
of finite
dimensional unital $*$-representations such that, for all  $k\geq 1$, 
$\oplus_{j=1}^{k}(\pi_j \oplus \hat{\pi}_j)$ is DPI. 
Let $\pi:=\oplus_{j \geq 1}\pi_j\oplus \hat{\pi}_j$ and $H:=\oplus_{j\geq 1}H_j\oplus \hat{H}_j$. 
To ease notation, for $k \geq 1$, let $\pi_{[k]}=\oplus_{j=1}^k\pi\oplus \hat{\pi}$.
Note that we have $\pi(A_1*A_2)\cap \mathbb{K}(H)=\{0\}$. Indeed, if $\pi(x)$ is compact then
$\lim_j \|(\pi_j\oplus\hat{\pi}_j )(x)\|=0$, since each representation is repeated infinitely
many times and we are considering a separating family we get $x=0$.

We will show that given any positive number  $\varepsilon$, there is a unitary $u$ on $\mathbb{U}(H)$ such 
that $\|u-\id_H\| < \varepsilon$ and  $\pi^{(1)}*(\Ad u \circ \pi^{(2)})$ is both irreducible and 
faithful. 
To do this, we will to construct a sequence $(u_k,\theta_k,F_k)_{k\geq 1}$ where:
\begin{enumerate}[(a)]
\item For all $k$, $u_k$  is a unitary in $\mathbb{U}(\oplus_{j=1}^kH_j\oplus \hat{H}_j)$ satisfying

\begin{eqnarray}\label{ConvergenceUnitaries}
\| u_k - \id_{\oplus_{j=1}^kH_j\oplus \hat{H}_j}\| < \frac{\varepsilon}{2^{k+1}}.
\end{eqnarray}

\item Letting 
\[
u_{(j,k)}=u_j \oplus \id_{H_{j+1} \oplus \hat{H}_{j+1}} \oplus \cdots  \oplus \id_{H_k\oplus \hat{H}_k}  
\]
and
\begin{equation}\label{Eqn:ExtendingU}
U_k=u_ku_{(k-1,k)} u_{(k-2,k)}\cdots u_{(1,k)}\,,
\end{equation}
the unital $*$-representation of $A_1*A_2 $ onto 
$\mathbb{B}\bigl( \oplus_{j=1}^k H_j\oplus \hat{H}_j \bigr)$, given by
\begin{eqnarray}\label{RelationWithPi}
\theta_k=  \pi_{[k]}^{(1)} * (\Ad U_k \circ \pi_{[k]}^{(2)} ),
\end{eqnarray}
is irreducible.

\item $F_k $ is a finite subset of the closed unit ball of $A_1*A_2$   and for all $y$ in the closed 
unit  ball of  $A_1* A_2$ there is an element $x$ in $F_k$ such that
\begin{eqnarray}\label{FiniteNet}
 \| \theta_k(x)-\theta_k(y) \| < \frac{1}{2^{k+1}}\,.
\end{eqnarray}

\item If $k\geq 2$, then for any element  $x$ in the union $\cup_{j=1}^{k-1} F_j$,  we have
\begin{eqnarray}\label{RememberingSetsF}
\|  \theta_{k}(x)-(\theta_{k-1} \oplus \pi_{k}\oplus \hat{\pi}_{k})(x) \|  < \frac{1}{2^{k+1}}\,.
\end{eqnarray}

\end{enumerate}
We construct such a sequence by recursion.

\noindent
{\em Step 1: Construction of $(u_1,\theta_1,F_1)$.}
Since $\pi\oplus \hat{\pi}$ is DPI, there is a unitary $u_1$ in $H_1\oplus \hat{H}_1$ 
such that  $\|u_1-\id_{H\oplus \hat{H}}\|< \frac{\varepsilon}{2^2}$ and 
$\pi_{[1]}^{(1)}*\Ad u_1 \circ \pi_{[1]}^{(2)}$ is
irreducible. Hence condition (\ref{ConvergenceUnitaries})  and (\ref{RelationWithPi})  trivially  hold.
Since $H_1 \oplus \hat{H}_1$ is finite dimensional,  
there is a finite set $F_1$ contained in the closed unit  ball of $A_1*A_2$ satisfying condition 
(\ref{FiniteNet}).  
At this stage there is no condition (\ref{RememberingSetsF}).

\noindent
{\em Step 2: Construction of $(u_{k+1},\theta_{k+1},F_{k+1})$ from $(u_j,\theta_j,F_j)$, $1\leq  j \leq k$.}
First, we are prove there exists  a unitary $u_{k+1}$ in $\mathbb{U}(\oplus_{j=1}^{k+1}H_j\oplus \hat{H}_j)$ 
such that $\| u_{k+1} - \id_{\oplus_{j=1}^{k+1}H_j \oplus \hat{H}_j} \| < \frac{\varepsilon}{2^{k+2}}$,
the unital  $*$-representation of $A_1*A_2$ into  
$\mathbb{B}\bigl( \oplus_{j=1}^{k+1} H_j \oplus \hat{H}_j \bigr)$ defined by 
\begin{equation}\label{DefinitioRepTheta}
\theta_{k+1}:= (\theta_k \oplus \pi_{k+1}\oplus \hat{\pi}_{k+1})^{(1)}*(\Ad u_{k+1} \circ 
(\theta_k \oplus \pi_{k+1}\oplus \hat{\pi}_{k+1} )^{(2)})  
\end{equation}
 is  irreducible and for any element $x$ in the union $\cup_{j=1}^k F_j$, the inequality  
$\|\theta_{k+1}(x)- (\theta_k \oplus \pi_{k+1} \oplus \hat{\pi}_{k+1})(x)   \|< \frac{1}{2^{k+1}}$, holds. 
By Remark~\ref{Remark:DPI}, $\theta_k \oplus \pi_{k+1} \oplus \hat{\pi}_{k+1}$ is 
DPI  so  Proposition~\ref{Perturbations&ApproximationOnGivenFiniteSets} 
assures the existence of such unitary  $u_{k+1}$.
Notice that, from construction,  conditions (\ref{ConvergenceUnitaries}) and  
(\ref{RememberingSetsF}) are satisfied.
A consequence of (\ref{RelationWithPi}) and 
(\ref{Eqn:ExtendingU})  is
\[
\theta_{k+1}=\pi_{[k+1]}^{(1)}*(\Ad U_{k+1} \circ \pi_{[k+1]}^{(2)}).
\]
Finite dimensionality of   $\oplus_{j=1}^{k+1}H_j\oplus\hat{H}_j $  guarantees the existence of a
finite set $F_{k+1}$ contained in the closed unit ball of $A_1*A_2$ satisfying condition  (\ref{FiniteNet}).
This completes Step 2.

Now consider the $*$-representations  
\begin{equation}\label{DefRepSigmak}
\sigma_k=\theta_k \oplus \bigoplus_{j\geq k+1} \pi_j \oplus \hat{\pi}_j.
\end{equation}
We now show there is a unital $*$-representation of  $\sigma:A_1*A_2\to\mathbb{B}(H)$, 
such that for  all $x$ in $A_1*A_2$, $\lim_k \|  \sigma_k(x) -\sigma(x)\|=0$.
If we extend the unitaries $u_k$ to all of $H$ via 
$\tilde{u}_k=u_k  \oplus_{j\geq k+1} \id_{H_j \oplus \hat{H}_j}$, 
then we  obtain
\begin{equation}
\sigma_k =  \pi^{(1)} *(\Ad \tilde{U}_k  \circ \pi^{(2)}) ,
\end{equation} 
where $\tilde{U}_k=\tilde{u}_k\cdots \tilde{u}_1$.
Thanks to condition (\ref{ConvergenceUnitaries}), we have
$$
\| \tilde{U}_k -\id_H \| \leq \sum_{j=1}^k\|  \tilde{u}_k - \id_H \| < \sum_{j=1}^k \frac{\varepsilon}{2^{k+1}},
$$
and for $l \geq 1$
$$
\| \tilde{U}_{k+l} - \tilde{U}_{k} \|  =\|  \tilde{u}_{k+l}\cdots \tilde{u}_{k+1}  - \id_H  \| 
\leq  \sum_{j=k+1}^{k+l} \frac{\varepsilon}{2^{j+1}}.
$$
Hence, Cauchy's criterion implies there is a unitary $u$ in $\mathbb{U}(H)$ such that the 
sequence $(\tilde{U}_k)_{k\geq 1}$ converges in norm to $u$ and $\|u-\id_H\|< \frac{\varepsilon}{2}$. Define
\begin{equation}\label{DefRepSigma}
\sigma=\pi^{(1)}*(\Ad u \circ \pi^{(2)}).
\end{equation}
From Proposition \ref{Perturbations&ApproximationOnGivenFiniteSets} we  have that for all $x$ in $A_1*A_2$,
\begin{equation}\label{PointWiseConvergenceSigma}
\lim_k \| \sigma_k(x)- \sigma(x)\| =0.
\end{equation}

Our next goal is to  show $\sigma$ is irreducible.  To ease notation let $A=A_1*A_2$.
We will show $\overline{\sigma(A)}^{SOT}=\mathbb{B}(H)$.
Take $T$ in $\mathbb{B}(H)$. With  no loss of generality we may assume 
$\|T\| \leq \frac{1}{2}$. Recall that a neighborhood basis for the
SOT topology around $T$ is given by the sets 
\[
\mathcal{N}_T(\xi_1,\dots, \xi_n;\varepsilon)=\{ S\in \mathbb{B}(H): \|S\xi_i -T\xi_i\|< \varepsilon
, i=1,\dots, n   \}
\]
where $\varepsilon >0$, $n \in \mathbb{N}$, and $\xi_1,\dots, \xi_n\in H$ are unit vectors.
We show that for any $\varepsilon>0$ and any unit vectors $\xi_1,\dots ,\xi_n$,
$\mathcal{N}_T(\xi_1,\dots, \xi_n;\varepsilon)\cap \sigma(A)$ is nonempty.
Let $P_k$ denote the orthogonal projection from $H$ onto
$\oplus_{j=1}^kH_j\oplus \hat{H}_j$. 
Take $k_1\geq 1$ such 
\[
\sum_{k \geq k_1}\frac{1}{2^k}< \frac{\varepsilon}{2^3}
\]
and for $k\geq k_1$, $1\leq i \leq n$,
\begin{eqnarray}
\|(\id_H-P_k)(\xi_i)\|< \frac{\varepsilon}{2^3} \label{Eqn:TailForxi}, \\
\|(\id_H - P_k)(T\xi_i)\| < \frac{\varepsilon}{2^3}.\label{Eqn:TailForT}
\end{eqnarray}
Since $P_k$ has finite rank and $\theta_k$ is irreducible,
there is $a$ in $A$, with $\|a\|\leq 1$ such that
\begin{equation}\label{Eqn:KadisonTransitivityAndTheta}
P_{k_1}TP_{k_1}(\xi_i)=\theta_{k_1}(a)(P_{k_1}(\xi_i))
\end{equation}
for $i=1,\dots ,n$.
We have
\begin{equation}\label{eq:thetasig}
\theta_{k_1}(a)(P_{k_1}(\xi_i))=\sigma_{k_1}(a)(P_{k_1}(\xi_i)).
\end{equation}
Take $x$ in $F_{k_1}$ such that
\begin{equation}\label{Eqn:RelationThetaAndFiniteNet}
\|\theta_{k_1}(a)- \theta_{k_1}(x)\| < \frac{1}{2^{k_1+1}}\,.
\end{equation}
We will show $\sigma(x)\in \mathcal{N}_T(\xi_1,\dots,\xi_n;\varepsilon )$.
To ease notation let $\xi_i=\xi$.
From (\ref{Eqn:TailForxi}), (\ref{Eqn:TailForT}), (\ref{Eqn:KadisonTransitivityAndTheta}) and~\eqref{eq:thetasig}, we deduce
\begin{eqnarray*}
\|T\xi - \sigma(x)\xi\| &\leq &\|T\xi - P_{k_1}TP_{k_1}\xi\| \\
 &+& \|P_{k_1}TP_{k_1}\xi - \sigma_{k_1}(a)\xi \|\\
&+& \| \sigma_{k_1}(a)\xi-\sigma(x)\xi\| \\
&<& \frac{3\varepsilon}{2^3} + \|\sigma_{k_1}(a)\xi - \sigma(x)\xi\|.
\end{eqnarray*}
For any $p\geq 1$ we have 
\begin{eqnarray*}
\sigma_{k_1}(a)\xi-\sigma(x)\xi&=& \sigma_{k_1}(a)\xi-\sigma_{k_1}(x)\xi \\
&+&\sum_{j=k_1}^{k_1+p}\big(\sigma_{j}(x)\xi -\sigma_{j+1}(x)\xi\big) \\
&+& \sigma_{k_1+p+1}(x)\xi- \sigma(x)\xi.
\end{eqnarray*}
Thus, from (\ref{Eqn:TailForxi}), \eqref{eq:thetasig}, (\ref{Eqn:RelationThetaAndFiniteNet}),
\eqref{DefRepSigmak} and (\ref{RememberingSetsF}) we deduce
\[
\|\sigma_{k_1}(a)\xi - \sigma(x)\xi\| < \frac{\varepsilon}{2} + \| \sigma_{k_1+p+1}(x)\xi - \sigma(x)\xi\| 
\]
hence
\[
\|\sigma_{k_1}(a)\xi - \sigma(x)\xi\| \leq  \frac{\varepsilon}{2}.
\]
We conclude $\sigma(x)$ lies in $\mathcal{N}_T(\xi_1,\dots, \xi_n;\varepsilon)$.

An application of Choi's technique  (see Theorem 6 in ~\cite{Choi}) will give us faithfulness of $\sigma$.  
Indeed, from construction, for all $x$ in $A$,
$\sigma(x)=\lim_k \sigma_k(x)$. Thus if each $\sigma_k$ is faithful then so is $\sigma$.
But faithfulness of $\sigma_k$ follows from the commutativity of the following diagram
$$
\xymatrix{A \ar[r]^{\pi} \ar[d]_{ \sigma_k } & \mathbb{B}(H) \ar[d]^{\pi_C} \\
\mathbb{B}(H) \ar[r]^-{\pi_C} & \mathbb{B}(H)/\mathbb{K}(H)
}
$$
(where $\pi_C$ denotes the quotient map onto the Calkin algebra), 
which in turn is implied by  (\ref{DefRepSigmak}).
\end{proof}

To obtain the following corollary, see Lemma~3.2 of~\cite{Bedos&Omland-Amenable}.

\begin{coro}
Assume $A_1$ and $A_2$ are nontrivial residually finite dimensional $C^*$-algebras with
$(\dim(A_1)-1)(\dim(A_2)-1) \geq 2$. Then $A_1*A_2$ is antiliminal and has an uncountable
family of pairwise inequivalent irreducible faithful $*$--representations. 
\end{coro}

We finish with a corollary derived from  Lemma 11.2.4 in ~\cite{Dixmier}.

\begin{coro}
Assume $A_1$ and $A_2$ are nontrivial residually finite dimensional $C^*$-algebras with
$(\dim(A_1)-1)(\dim(A_2)-1) \geq 2$. Then pure states of  $A_1*A_2$ are w*-dense in the state space.
\end{coro}


\begin{bibdiv}
\begin{biblist}

\bib{Bedos&Omland-Amenable}{article}{
   author={B{\'e}dos, Erik},
   author={Omland, Tron {\AA}.},
   title={Primitivity of some full group $C^\ast$-algebras},
   journal={Banach J. Math. Anal.},
   volume={5},
   date={2011},
   pages={44--58},
}

\bib{Bedos&Omland-Modular}{article}{
   author={B{\'e}dos, Erik},
   author={Omland, Tron {\AA}.},
   title={The full group C$^*$--algebra of the modular group is primitive},
   journal={Proc. Amer. Math. Soc.},
   volume={140},
   date={2012},
   pages={1403-1411},
}

\bib{Blackadar}{book}{
   author={Blackadar, B.},
   title={Operator algebras},
   series={Encyclopaedia of Mathematical Sciences},
   volume={122},
   note={Theory of $C^*$-algebras and von Neumann algebras;
   Operator Algebras and Non-commutative Geometry, III},
   publisher={Springer-Verlag},
   place={Berlin},
   date={2006},
}

\bib{Choi}{article}{
   author={Choi, Man Duen},
   title={The full $C^{\ast} $-algebra of the free group on two generators},
   journal={Pacific J. Math.},
   volume={87},
   date={1980},
   pages={41--48},
}

\bib{Davidson}{book}{
   author={Davidson, Kenneth R.},
   title={$C^*$--algebras by example},
   series={Fields Institute Monographs},
   volume={6},
   publisher={American Mathematical Society},
   place={Providence, RI},
   date={1996},
}

\bib{Dixmier}{book}{
  author={Dixmier, Jacques},
  title={Les C$^*$--alg\`ebres et leurs repr\'esentations},
  publisher={Gauthiers--Villars},
  year={1969}
}

\bib{Exel&Loring}{article}{
   author={Exel, Ruy},
   author={Loring, Terry A.},
   title={Finite-dimensional representations of free product $C^*$--algebras},
   journal={Internat. J. Math.},
   volume={3},
   date={1992},
   pages={469--476},
}

\bib{Helgason}{book}{
   author={Helgason, Sigurdur},
   title={Differential geometry, Lie groups, and symmetric spaces},
   series={Pure and Applied Mathematics},
   volume={80},
   publisher={Academic Press Inc. [Harcourt Brace Jovanovich Publishers]},
   place={New York},
   date={1978},
}

\bib{Kirillov}{book}{
   author={Kirillov, Alexander, Jr.},
   title={An introduction to Lie groups and Lie algebras},
   series={Cambridge Studies in Advanced Mathematics},
   volume={113},
   publisher={Cambridge University Press},
   place={Cambridge},
   date={2008},
}

\bib{Murphy}{article}{
   author={Murphy, G. J.},
   title={Primitivity conditions for full group $C^\ast$-algebras},
   journal={Bull. London Math. Soc.},
   volume={35},
   date={2003},
   pages={697--705},
}


\bib{Omland-PrimitivityConditionsTwistedGroups}{article}{
   author={Omland, Tron {\AA}.},
   title={Primeness and primitivity conditions for twisted group C$^*$--algebras },
   eprint={http://arxiv.org/abs/1204.4259v1}
}

\bib{Pedersen}{article}{
   author={Pedersen, Gert Kjaerg{\.a}rd},
   title={Measure theory for $C^{\ast} $ algebras. II},
   journal={Math. Scand.},
   volume={22},
   date={1968},
   pages={63--74},
}

\bib{Raeburn&Sinclair}{article}{
   author={Raeburn, Iain},
   author={Sinclair, Allan M.},
   title={The $C^*$-algebra generated by two projections},
   journal={Math. Scand.},
   volume={65},
   date={1989},
   pages={278--290},
}

\bib{Takesaki}{book}{
   author={Takesaki, M.},
   title={Theory of operator algebras. I},
   series={Encyclopaedia of Mathematical Sciences},
   volume={124},
   note={Reprint of the first (1979) edition;
   Operator Algebras and Non-commutative Geometry, 5},
   publisher={Springer-Verlag},
   place={Berlin},
   date={2002},
}

\bib{Yoshizawa}{article}{
   author={Yoshizawa, Hisaaki},
   title={Some remarks on unitary representations of the free group},
   journal={Osaka Math. J.},
   volume={3},
   date={1951},
   pages={55--63},
}

\end{biblist}
\end{bibdiv}
\end{document}